\newif\ifsingle

\ifsingle
\documentclass[11pt,draftclsnofoot, onecolumn]{IEEEtran}		
\else		
\documentclass[10pt,journal,compsoc]{IEEEtran}
\fi

\usepackage{times}
\usepackage{subcaption}
\usepackage{amsmath,dsfont}
\usepackage{amssymb,amsthm}
\usepackage{epsfig,verbatim} 
\usepackage{setspace}
\usepackage{color}
\usepackage{cite}
\usepackage{epstopdf}
\usepackage{graphics}
\usepackage{accents}
\usepackage{acronym}
\usepackage[bookmarks,colorlinks]{hyperref}
\usepackage{booktabs}
\usepackage{mathtools}
\usepackage{enumitem}
 \usepackage{diagbox}

\usepackage[ruled,linesnumbered,vlined]{algorithm2e}
\SetKwInput{KwData}{\textbf{Init}}

\definecolor{NewColor}{rgb}{0,0,0}

\newcommand{\myVec}[1]{{\boldsymbol{#1}}}
\newcommand{\myMat}[1]{{\boldsymbol{#1}}}
\newcommand{\mySet}[1]{\mathcal{#1}}

\acrodef{mse}[MSE]{mean-squared error}
\acrodef{iot}[IOT]{Internet of Things}
\acrodef{rmse}[RMSE]{root mean squared error}
\acrodef{rmspe}[RMSPE]{root mean squared periodic error}
\acrodef{mmse}[MMSE]{{minimum mean-squared error}}
\acrodef{lmmse}[LMMSE]{{linear} MMSE}
\acrodef{mle}[MLE]{maximum likelihood estimation}
\acrodef{snr}[SNR]{signal-to-noise ratio}
\acrodef{cnn}[CNN]{convolutional neural network}
\acrodef{rnn}[RNN]{recurrent neural network}
\acrodef{dnn}[DNN]{deep neural network}
\acrodef{gnn}[GNN]{graph neural network}
\acrodef{admm}[ADMM]{alternating direction method of multipliers}
\acrodef{dadmm}[D-ADMM]{distributed alternating direction method of multipliers}
\acrodef{sgd}[SGD]{stochastic gradient descent}
\acrodef{dlasso}[D-LASSO]{distributed least absolute shrinkage and selection operator}
\acrodef{dlr}[D-LR]{distributed linear regression}
\acrodef{fl}[FL]{federated learning}

\IEEEoverridecommandlockouts
\title{Distributed Learn-to-Optimize: 
Limited Communications Optimization over Networks via Deep Unfolded Distributed ADMM}

\author{
\IEEEauthorblockN{Yoav Noah and Nir Shlezinger
\thanks{
Parts of this work were presented in the 2023 IEEE International Conference on Acoustics Speech, and Signal Processing (ICASSP) as the paper \cite{noah2023distributed}. 
The authors are with the School of ECE, Ben-Gurion University of the Negev, Israel (e-mail: yoavno@post.bgu.ac.il; nirshl@bgu.ac.il). 
}}}

\begin{document}

\maketitle
\pagestyle{plain}
\thispagestyle{plain}

%
%
\begin{abstract} 
Distributed optimization is a fundamental framework for collaborative inference over networks. The operation is modeled as the joint minimization of a shared objective which typically depends on local observations. Distributed optimization algorithms, such as the  \ac{dadmm}, iteratively combine local computations and message exchanges. A main challenge associated with distributed optimization, and particularly with \ac{dadmm}, is that it requires a  large number of communications to reach consensus. 
In this work we propose {\em unfolded \ac{dadmm}}, which follows the emerging deep unfolding methodology to enable \ac{dadmm} to operate reliably with a predefined and small number of messages exchanged by each agent. Unfolded \ac{dadmm} fully preserves the operation of \ac{dadmm}, while leveraging data to tune the hyperparameters of each iteration. These hyperparameters can either be agent-specific, aiming at achieving the best performance within a fixed number of iterations over a given network, or shared among the agents, allowing to learn to distributedly optimize over different networks. 
We specialize unfolded \ac{dadmm} for two representative settings: a distributed sparse recovery setup, and a distributed machine learning learning scenario.
Our numerical results demonstrate that the proposed approach dramatically reduces the number of communications utilized by \ac{dadmm}, without compromising on its performance.  
\end{abstract}
\acresetall

\section{Introduction}
\label{sec:intro}
The proliferation of sophisticated electronic devices results in data, e.g., locally sensed signals, being divided among many different agents. By collaborating with each other, the users can jointly extract desired information or learn a machine learning models from the divided data in a decentralized fashion. This is achieved by {\em distributed optimization}~\cite{nedic2020distributed}, leveraging the communications capabilities of the agents, without relying on centralized data pooling implemented by traditional centralized systems.
A shared objective is iteratively optimized in a distributed manner by having each agent process the data locally to solve a local optimization problem, and then update its neighbors, repeating until consensus is achieved~\cite{xin2020general}. Various different algorithms implement distributed optimization, see survey in \cite{yang2019survey}.  A popular and common distributed optimization algorithms is  \ac{dadmm}~\cite{boyd2011distributed,mota2013d}, which can often be guaranteed to converge into optimal consensus~\cite{ han2012note,shi2014linear,makhdoumi2017convergence}. 

Distributed optimization algorithms, including \ac{dadmm}, typically require multiple iterations to reach consensus. Each iteration involves not only local computations, but also message exchanges between the participating agents.  
The latter implies excessive communications, possibly inducing notable delays on the overall optimization procedure, limiting its real-time applicability, and can be costly in power and spectral resources \cite{akyildiz2002wireless, fischione2011design}.
This raises a notable challenge that is not encountered in centralized optimization, and can be a limiting factor in multi-agent systems that are inherently limited in power and communications, such as \ac{iot} and sensor networks~\cite{rabbat2004distributed, shlezinger2022collaborative}.  

Two leading approaches are proposed in the literature to carry out distributed optimization with reduced communications. The first approach aims at deriving alternative optimization algorithms with improved convergence rates~\cite{jakovetic2014fast, chang2014multi, mokhtari2016decentralized,aybat2017distributed, lan2020communication}. For instance, the primal dual method of multipliers~\cite{zhang2017distributed} is closely related to \ac{dadmm}, and can be shown to converge more rapidly (assuming it converges)~\cite{sherson2018derivation}. Alternatively, one can reduce the communication overhead by quantizing the exchanged messages, such that each message is constrained to a fixed and small number of bits \cite{aysal2008distributed, kar2009distributed, zhu2015quantized,shlezinger20221quantized,cohen2021serial}. 
However, these conventional approaches are typically studied in the context of convergence rate, which is an asymptotic property. In practice, one is often interested in operating with a fixed latency, and there is a concrete need to enable reliable low-latency distributed optimization techniques to be  carried out with a small and fixed number of communication rounds. 

We are witnessing a growing interest in machine learning, and particularly in deep learning, beyond its traditional domains of computer vision and natural language processing. 
Deep learning systems employ data-driven parameterized \acp{dnn} that learn their mapping from data~\cite{Goodfellow_Deep_Learning}.
\textcolor{NewColor}{
In the context of networks and multi-agent systems, deep learning tools were shown to facilitate various challenging tasks, such resource allocation~\cite{samikwa2022ares,lee2022radio}, anomaly detection~\cite{alani2022phishnot,cohen2024sinr}, and network management~\cite{zhang2023ai,khoshkholghi2022edge}, see also~\cite{cerquitelli2023machine}. For tasks that can be modeled as distributed inference or decision making, existing literature on the integration of \acp{dnn} typically focuses on centralized scenarios~\cite{song2022networking,chen2024learned}. There, the agents communicate in a star topology directly with a cloud server~\cite{lee2023task} or with a specific node that orchestrates the procedure~\cite{malka2022decentralized}. This is as opposed to distributed optimization such as \ac{dadmm}, which is concerned with agents communicating over arbitrary mesh topologies and aim to  reach consensus in a decentralized manner.}

\textcolor{NewColor}{As opposed to tasks that can be represented as centralized or remote inference, the application of deep learning for distributed optimization is far less explored. 
Existing families of \ac{dnn} architectures vary based on the type of data, with traditional fully connected networks designed to process unstructured vectors of fixed size, while, e.g., convolutional and recurrent neural networks are suitable for images and time sequences, respectively~\cite[Ch. 9-10]{Goodfellow_Deep_Learning}. The local observations based on which distributed optimization is carried out can be represented as graph signals~\cite{vatter2023evolution}, with the leading \ac{dnn} architecture for processing such signals being  \acp{gnn}~\cite{wu2020comprehensive}. The learned operation of \acp{gnn} involves iterative message exchange between agents and local processing, where each \ac{gnn} layer implements a single message exchange round.} Accordingly, various \ac{gnn} architectures were recently employed for tasks involving distributed optimization~\cite{shen2020graph, rusek2020routenet} and learning~\cite{hadou2023stochastic}. 
\textcolor{NewColor}{
Alternative designs equip each agent with a different \ac{dnn} whose dimensions depend on the number of agents~\cite{lee2019deep}, at the cost of limited scalability and robustness to topology variations.}
Although such \acp{dnn} can learn to implement distributed optimization with a fixed and limited number of communication rounds, they lack the interpretability of conventional distributed optimization, and often require large volumes of data for learning purposes.

While deep learning methods are traditionally considered to replace principled inference based on mathematical modelling, they can also be integrated with classical inference~\cite{shlezinger2020model,chen2022learning,shlezinger2023model}. In particular, deep learning tools can be employed to {\em learn-to-optimize}~\cite{chen2022learning}, i.e., improve iterative optimization in performance and run-time~\cite{shlezinger2022model}. This can be realized using the deep unfolding methodology~\cite{monga2021algorithm}, that leverages data to optimize the performance of an iterative inference rule with a fixed number of iterations. Recent works have considered the distributed optimization of centralized deep unfolded mappings~\cite{mogilipalepu2021federated} and the integration of deep unfolding for centralized aggregation in distributed learning~\cite{nakai2022deep}. Yet, deep unfolding is currently considered in the context of centralized optimization, where its gains are mostly in run-time and in possibly increased abstractness of the resulting inference rule~\cite{shlezinger2022model}, motivating its exploration for facilitating distributed optimization with a fixed and small number of communication rounds. 
In this work we explore the usage of deep unfolding to facilitate distributed optimization with a fixed and limited communication budget. We focus on the \ac{dadmm} algorithm, being a popular and common distributed optimization algorithm, aiming to show the usefulness of the combination of deep unfolding with distributed optimization. To that aim,  we unfold \ac{dadmm} to operate with a predefined  fixed number of communication rounds. Then, we leverage data to tune the optimization hyperparameters, converting \ac{dadmm} into a trainable discriminative machine learning model~\cite{shlezinger2022discriminative} such that the resulting distributed mapping fits the training data. Our proposed {\em unfolded \ac{dadmm}} thus learns from data to achieve improved performance within the fixed number of iterations.  

We consider two parameterizations for the unfolded architecture. The first learns {\em agent-specific hyperparameters}, i.e., allowing them to vary not only between iterations, but also between agents. The second learns {\em shared hyperparameters}, where all agents use the same hyperparameters that vary between iterations. As each setting learned with shared hyperparameters is a special case of the agent specific case, this latter allows achieving improved performance.  However, learning shared hyperparameters is not limited to a given set of agents, and thus a mapping learned for a given network can generalize to different networks. In both cases, the resulting unfolded \ac{dadmm} thus completely preserves the  operation of the conventional \ac{dadmm}, while allowing it to achieve improved performance within a predefined and small communication budget. 

We showcase the application of unfolded \ac{dadmm} in two representative case studies: $(i)$ distributed sparse recovery; and $(ii)$ the distributed learning of a linear regression model. For each setting, we specialize the formulation of \ac{dadmm} and its unfolded architecture. Our experimental studies show that the application of deep unfolding yields notable improvements in both case studies. In particular, we demonstrate  reductions by factors varying between $\times 8$ and up to $\times 154$ in the amount of communications compared with the conventional \ac{dadmm}, without degrading its performance and while completely preserving its interpretable operation. We also show that by following a principled optimization algorithm, our unfolded \ac{dadmm} notably outperforms \ac{gnn} architectures trained for the same task. 


The rest of this work is organized as follows: Section~\ref{sec:Model} formulates the generic distributed optimization setup and recalls the conventional \ac{dadmm}. The proposed unfolded D-ADMM algorithm is derived in Section~\ref{sec:Unfolded}. The considered case studies of distributed sparse recovery and distributed linear regression learning are reported in Sections~\ref{sec:LASSO-Model}-\ref{sec:LinearRegModel}, respectively, along with their corresponding experimental study. Finally,  Section~\ref{ssec:conclusions} provides concluding remarks.

%
\section{System Model}\label{sec:Model}
%
This section provides the necessary background needed to derive our proposed unfolded \ac{dadmm} in Section~\ref{sec:Unfolded}. In particular, we first formulate the generic distributed optimization problem in Subsection~\ref{ssec:problem}, after which we recall \ac{dadmm} in Subsection~\ref{ssec:d-admm}.

%
\subsection{Generic Distributed Optimization Problem Formulation}\label{ssec:problem} 
%
We consider a set of $P$ agents indexed in $1,\ldots,P$. Each agent of index $p \in 1,\ldots,P$ has access to a local data  $\myVec{b}_p \in \mathbb{R}^m$, representing, e.g., locally acquired observations. Based on the local data, the agents aim at jointly solving an optimization problem of the form 
\vspace{-5pt}
\begin{equation}
    \label{eqn:opt-problem} 
    \mathop{\arg \min}\limits_{\myVec{\Bar{y}} 
    } \quad  \sum_{p=1}^{P} f_p(\myVec{\Bar{y}}; \myVec{b}_p).
\end{equation}
In \eqref{eqn:opt-problem}, the vector $\myVec{\Bar{y}} \in \mathbb{R}^{n}$ is the optimized variable, and $f_p: \mathbb{R}^{n} \times \mathbb{R}^m \mapsto \mathbb{R}^+$ is the $p$th objective. 

The generic formulation in \eqref{eqn:opt-problem}  can be viewed as a centralized optimization problem with a decomposable objective. Yet, the distributed nature of \eqref{eqn:opt-problem} stems from the fact that each observation $\myVec{b}_p$ is available solely to the $p$th agent, and there is no data pooling, i.e.,  the gathering of $\{\myVec{b}_p\}_{p=1}^{P}$ by a centralized data fusion entity.

%
%

Each agent can communicate with its neighbours. The resulting communication network is modeled as an undirected connected graph  $G(\mySet{V}, \mySet{E})$, where $\mySet{V} = \{1, \ldots, P\}$ is the set of nodes and $\mySet{E} \subseteq \mySet{V} \times \mySet{V}$ is the set of edges. 
\textcolor{NewColor}{This modelling follows the common representation of multi-agent networks in distributed optimization, see, e.g.,~\cite{boyd2011distributed,mota2013d,yang2019survey,shi2014linear,makhdoumi2017convergence}.}
The graph is proper colored, i.e., there exist non-overlapping sets $\{\mySet{V}_c\}_{c=1}^C$ where $\cup_{c=1}^C\mySet{V}_c = \mySet{V}$ and no two vertices in each $\mySet{V}_c$ share an edge. This definition specializes uncolored graphs for $C=1$. 

Our goal in this paper is to tackle the optimization problem formulated in \eqref{eqn:dadmmopt-problem} under a  scenario with limited predefined communication and maximal latency. Here, each node is allowed to transmit up to $T$ messages to each of its neighbours. Furthermore, to enable deriving rapid communications-limited optimizers, we also have access to a labeled data set comprised of $L$ pairs of local data sets along the corresponding desired optimization variable. The data set is given by
\begin{equation}
    \label{eqn:dataset}
    \mySet{D} = \left\{\left(\left\{\myVec{b}_{l,p}\right\}_{p=1}^{P}, \myVec{\Bar{y}}_{l}\right)\right\}_{l=1}^{L},
\end{equation}
where for each sample of index $l$, $\myVec{b}_{l,p} \in \mathbb{R}^{m}$ is the observations vector at the $p$th node, and $\myVec{\Bar{y}}_{l} \in \mathbb{R}^{n}$ is the desired optimization outcome.  

\subsection{D-ADMM Algorithm}\label{ssec:d-admm} 

 A popular method for tackling distributed optimization problems as formulated in \eqref{eqn:opt-problem} is the \ac{dadmm} algorithm~\cite{mota2013d}. To formulate \ac{dadmm}, we first cast \eqref{eqn:opt-problem} as a consensus setup. 
 This is done by assigning copies of the global parameter $\myVec{\bar{y}}$ to each node while constraining all copies to be equal. Specifically, since the network is connected, the constraint is equivalently  imposed for each pair of neighbouring agents~\cite{zhu2009distributed}. The resulting distributed optimization problem is reformulated as
 \vspace{-0.1cm}
\begin{align}
 \label{eqn:dadmmopt-problem}
        &\mathop{\arg \min}\limits_{\myVec{y}_{1}, \ldots, \myVec{y}_{P}
    } \quad  \sum_{p=1}^{P} f_{p}(\myVec{y}_{p}; \myVec{b}_{p}) \\
    &\text{subject to}  \quad \myVec{y}_{p} = \myVec{y}_{j}, \forall j \in \mySet{N}_p.\notag
 \vspace{-0.1cm}
\end{align}
In \eqref{eqn:dadmmopt-problem},  $\mySet{N}_{p}$ denotes the set of all the neighbors of node $p$ in $G(\mySet{V}, \mySet{E})$, and $\myVec{\bar{y}} = (\myVec{y}_{1}, \dots, \myVec{y}_{p}) \in (\mathbb{R}^{n})^{P}$ is the optimization variable. 
The equivalent formulation in \eqref{eqn:dadmmopt-problem} represents \eqref{eqn:opt-problem} as a set of $P$ local optimziation problems, that are constrained to have identical solutions among the different agents. 


\ac{dadmm} tackles the constrained optimization in \eqref{eqn:dadmmopt-problem} by formulating the augmented  Lagrangian for each agent. For the $p$th agent, the augmented Lagrangian is given by 
\begin{align}
    &\mathcal{L}_p(\myVec{y}_{p}, \{\myVec{y}_{j}\}_{j\in\mySet{N}_p}, \myVec{\lambda}_p) \triangleq 
    f_{p}(\myVec{y}_{p};\myVec{b}_{p})
    \notag \\
    &\qquad \qquad + \sum_{j \in \mySet{N}_p} \myVec{\lambda}_p^T(\myVec{y}_{p} - \myVec{y}_{j}) + \frac{\rho}{2}\|\myVec{y}_{p} - \myVec{y}_{j}\|_2^2,
    \label{eqn:AugLag}
\end{align}
where $\rho>0$ is a fixed hyperparameter and $\myVec{\lambda}_p$ is the dual variable of node $p$.
\ac{dadmm} then has each agent alternate between $(i)$ minimizing \eqref{eqn:AugLag} with respect to the local optimization variable $\myVec{y}_{p}$; $(ii)$ share this update with its neighbours; and $(iii)$ maximize  \eqref{eqn:AugLag} with respect to the local dual variable $\myVec{\lambda}_p$.

The implementation of the minimization and maximization steps typically depends on the specific objective functions. Following~\cite{mota2013d}, we focus on the generic implementation of these operations by gradient descent and gradient ascent steps, respectively. By letting $\{\myVec{y}_{j,p}\}_{j\in\mySet{N}_p}$ denote the copies of $\{\myVec{y}_{j}\}$ available to the $p$th agent, the $k$th iteration computes
\begin{align}
    \myVec{y}_{p}^{(k+1)} &= \myVec{y}_{p}^{(k)} -\alpha \nabla_{\myVec{y}_{p}} \mathcal{L}_p\left(\myVec{y}_{p}^{(k)}, \{\myVec{y}_{j,p}\}_{j\in\mySet{N}_p}, \myVec{\lambda}_p^{(k)}\right),
    \label{eqn:PrimalStep}
\end{align}
by each agent $p$, with $\alpha>0$ being a step-size.
The agent then shares $\myVec{y}_{p}^{(k+1)}$ with its neighbours, who update their local copies. The iteration is concluded by having all agents update their dual variable with step-size $\eta>0$ via
\begin{align}
    \myVec{\lambda}_{p}^{(k+1)} &= \myVec{\lambda}_{p}^{(k)} +\eta \nabla_{\myVec{\lambda}_{p}} \mathcal{L}_p\left(\myVec{y}_{p}^{(k+1)}, \{\myVec{y}_{j,p}\}_{j\in\mySet{N}_p}, \myVec{\lambda}_p^{(k)}\right). 
    \label{eqn:DualStep}
\end{align}

The fact that the graph is colored allows some of the computations in the above iterations to be carried out in parallel.  The resulting procedure is summarized as Algorithm~\ref{alg:Algo1}, where  all nodes in the same color group  simultaneously update \eqref{eqn:PrimalStep}, and the operation is repeated until a stopping criteria is reached, e.g., convergence is achieved or a maximal number of iterations is exhausted. While this reduces the run-time of each iteration, it does not reduce the amount of message exchanges, which is dictated by the number of iterations required to reach a the stopping criterion, and is often very large (in the order of hundreds or thousands of iterations to converge).

\begin{algorithm}  
	\caption{D-ADMM}
	\label{alg:Algo1}
	\KwData{$\forall p \in \mySet{V}$, set $\myVec{\lambda}_{p}^{(1)}$, $\myVec{y}_{p}^{(1)}$,  $\{\myVec{y}_{j,p}\}$ to zero;  $k=1$}
	\Repeat{stopping criteria is reached}{
	    \For{$c=1,2,\ldots,C$}{
		\For{all $p \in \mySet{V}_{c}$ [in parallel]}{
		        Update 
		        $\myVec{y}_{p}^{(k+1)}$ via \eqref{eqn:PrimalStep}; \\ 
		        \textrm{Send} $\myVec{y}_{p}^{(k+1)}$ to $\mySet{N}_{p}$;
		        } 
	    }
	    \For{all $p \in \mySet{V}$ [in parallel]}{
	         Update 
		        $\myVec{\lambda}_{p}^{(k+1)}$ via \eqref{eqn:DualStep};
	        }
	   $k \xleftarrow{} k+1$
	    } 
\end{algorithm}

The above formulation of \ac{dadmm} is given for a generic problem, i.e., with explicitly stating the local objective functions $\{f_p(\cdot)\}$. Two special cases of the application of \ac{dadmm} for distributed sparse recovery and for distributed learning of a linear regression model are detailed in Sections~\ref{sec:LASSO-Model} and~\ref{sec:LinearRegModel}, respectively.

\section{Deep Unfolded D-ADMM}
\label{sec:Unfolded}

Although \ac{dadmm} is a suitable algorithm for tackling \eqref{eqn:opt-problem} with the objective \eqref{eqn:AugLag}, its direct application is likely to yield inaccurate estimates when applied with a fixed and small communication budget $T$, i.e., when allowed to run up to $T$ iterations. In this section we describe how we leverage the available data set in \eqref{eqn:dataset} to enable \ac{dadmm} to operate reliably and rapidly via deep unfolding. As deep unfolding converts an iterative algorithm into a sequential  discriminative model~\cite{shlezinger2022discriminative}, we begin by describing the  trainable architecture in Subsection~\ref{ssec:Architecture}. Then, we detail the training procedure in Subsection~\ref{ssec:Training}, after which we provide a discussion in Subsection~\ref{ssec:discussion}. 


%
\subsection{Trainable Architecture}
\label{ssec:Architecture}
We unfold \ac{dadmm} by fixing its number of iterations to be $T$, thus meeting the communication budget imposed in Section~\ref{sec:Model}. In order to enable the resulting algorithm to infer reliably, we note that \ac{dadmm} tackles the objective in~\eqref{eqn:dadmmopt-problem} by introducing three hyperparameters: the regularization coefficient $\rho$, and the primal-dual step-sizes $\alpha$ and $\eta$. While the exact setting of $(\rho, \alpha, \eta)$  typically does not affect the algorithm outcome when allowed to run until convergence (under some conditions~\cite{boyd2011distributed,han2012note}),  they have notable effect when using a fixed number of iterations~\cite{shlezinger2022model}. 

Our proposed unfolded \ac{dadmm} builds upon the insight that a proper setting of $(\rho, \alpha, \eta)$ can notably improve the performance of \ac{dadmm} when constrained to use $T$ iterations. We thus treat \ac{dadmm} with $T$ iterations as  discriminative machine learning model by considering two possible settings of trainable parameters: {\em agent-specific hyperparemeters} and {\em shared hyperparameters}.

\subsubsection{Agent-Specific Hyperparameters}\label{ssec:diff-hyper} 
Here, we allow $(\rho, \alpha, \eta)$ to {\em vary between agents} and between iterations, while treating them as trainable parameters. By doing so, each iteration of \ac{dadmm} can be viewed as a layer of a $T$-layered trainable model; the parameters of the $k$th layer are given by 
\begin{equation}
\label{eqn:agentSepc}
 \myVec{\theta}_{k} = \{\rho_{p}^{(k)}, \alpha_{p}^{(k)}, \eta_{p}^{(k)}\}_{p=1}^{P},   
\end{equation}
where ${\rho}_{p}^{(k)}$, ${\alpha}_{p}^{(k)}$, and ${\eta}_{p}^{(k)}$ are used as the setting of $\rho$, $\alpha$ and $\eta$, by the $p$th agent at the $k$th iteration of Algorithm~\ref{alg:Algo1}. 

The agent-specific approach parameterization allows each agent to use different hyperparameters that also change between iterations. As such, it provides increased flexibility to the operation of \ac{dadmm} with $T$ iterations for a given graph $\mySet{G}(\mySet{V},\mySet{E})$, which is exploited by learning these hyperparameters from data (see Subsection~\ref{ssec:Training}). However, the fact that this form of unfolded \ac{dadmm} assigns a specific set of hyperparameters to each individual agent implies that it should be trained and applied on the same set of agents. Furthermore, the number of trainable parameters grows with the number of agents, which make training more complicated and computationally complex, particularly when dealing with massive networks. 

\subsubsection{Shared Hyperparameters}\label{ssec:same-hyper} 
In order to decouple the parameterization of the unfolded \ac{dadmm} from the graph $\mySet{G}(\mySet{V},\mySet{E})$, we also consider a setting with shared hyperparameters. Here, the hyperparameters vary between iterations but are {\em shared between agents}. Accordingly, each \ac{dadmm} iteration is viewed as a layer of a $T$-layered trainable model, 
 where the parameters of the $k$th layer are given by 
\begin{equation}
\label{eqn:shared}
 \myVec{\theta}_{k} = \{\rho^{(k)}, \alpha^{(k)}, \eta^{(k)}\},   
\end{equation}
with ${\rho}^{(k)}$, ${\alpha}^{(k)}$, and ${\eta}^{(k)}$ being the setting of $\rho$, $\alpha$ and $\eta$ employed at the $k$th iteration of Algorithm~\ref{alg:Algo1}.

The parameterization in \eqref{eqn:shared} is clearly a constrained case of that allowed in \eqref{eqn:agentSepc}, and thus every configuration learned with shared hyperparaemters can also be learned with agent specific hyperparameters. This indicates that this conversion of \ac{dadmm} into a trainable machine learning model is not expected to achieve improved performance compared with the agent-specific setting, when training and inference are carried out on the same graph. However, the main benefit of the shared hyperparameters approach is that it is not sensitive to the graph specificities. As such, the same learned hyperparameters can be employed on different graphs, as we numerically demonstrate in Section~\ref{sec:LASSO-Model}. Furthermore, it allows our unfolded \ac{dadmm} can be trained end-to-end on a small graph and then have trained hyperparametres used on a bigger graph, an approach which is known to facilitate learning on large graphs in the context of \acp{gnn} \cite{cervino2023training}. 



%
\vspace{-0.2cm}
\subsection{Training Procedure}
\label{ssec:Training}
 
The unfolding of \ac{dadmm} into a $T$-layered model yields an architecture with trainable parameters $\myVec{\theta} = \{\myVec{\theta}_k\}_{k=1}^{T}$. In particular,  $\myVec{\theta}$ is comprised $3 \cdot P \cdot T$  trainable parameters for agent-specific hyperparameters or of $3 \cdot T$ parameters for shared hyperparameters (though one can possibly introduce additional trainable parameters arising, e.g., from the objective, as we do in Section~\ref{sec:LASSO-Model}). 
We use the data set $\mySet{D}$ in \eqref{eqn:dataset} to tune $\myVec{\theta}$ such that the local decisions obtained by the agents, i.e., $\{\bar{\myVec{y}}_p\}$, accurately match the desired consensus output $\bar{\myVec{y}}$. 

In particular, we train unfolded \ac{dadmm} using the \ac{mse} loss. Let $\myVec{y}^{(k)}_{p}\big(\{\myVec{b}_i\}_{i=1}^P;\myVec{\theta}\big)$ be the estimate produced at the $k$th iteration at node $p$ when applying \ac{dadmm} with parameters $\myVec{\theta}$ and observations $\{\myVec{b}_i\}_{i=1}^P$. The resulting training loss is  %
\begin{equation}
    \label{eqn:mseloss}
    \begin{aligned}
        \mathcal{L}_{\mySet{D}}(\myVec{\theta})\! =\! \frac{1}{|\mySet{D}| P} \! \sum_{(\{\myVec{b}_{l,i}\},\myVec{\Bar{y}}_{l})\in \mySet{D} }\sum_{p=1}^{P} \lVert \myVec{y}^{(T)}_{p}\!\big(\{\myVec{b}_{l,i}\};\myVec{\theta}\big)\! -\!\myVec{\Bar{y}}_{l} \rVert_2^{2}.
    \end{aligned}
\end{equation}


We optimize the parameters $\myVec{\theta}$ based on  \eqref{eqn:mseloss}
using deep learning optimization techniques based on, e.g., mini-batch stochastic gradient descent. 

When training on a large number of agents it is likely to suffer from hardware limitations, hence the training process is done using a sequential training method~\cite{shlezinger2020deepsic}. Unlike conventional batch training, which requires processing large volumes of data all at once, sequential training takes a more incremental approach. It divides the training process into smaller segments, allowing the model to be updated progressively. By breaking down the training process into $t$-layers segments, where $t < T$, we mitigate the memory and computation bottlenecks, reduce the immediate strain on the hardware. The $t$th layer output is a soft estimate of the \ac{mse} loss for each $t \in T$. Therefore, each building block of unfolded \ac{dadmm} can be trained individually, by minimizing the \ac{mse} loss. To formulate this objective, let $\myVec{y}^{(t)}_{p}\big(\{\myVec{b}_i\}_{i=1}^P;\myVec{\theta}_{t}\big)$ be the estimate produced at the $t$th iteration at node $p$ when applying \ac{dadmm} with parameters $\myVec{\theta}_{t}$ and observations $\{\myVec{b}_i\}_{i=1}^P$, where $\myVec{\theta}_{t}$ represents the hyperparameters at the $t$th layer. The resulting training loss is
\begin{equation}
    \label{eqn:mseloss-sequntial}
    \begin{aligned}
        \!\!\!\!\mathcal{L}_{\mySet{D}}(\myVec{\theta}_{t})\! =\! \frac{1}{|\mySet{D}| P} \! \sum_{(\{\myVec{b}_{l,i}\},\myVec{\Bar{y}}_{l})\in \mySet{D} }\sum_{p=1}^{P} \lVert \myVec{y}^{(t)}_{p}\!\big(\{\myVec{b}_{l,i}\};\myVec{\theta}_{t}\big)\! -\!\myVec{\Bar{y}}_{l} \rVert_2^{2}.
    \end{aligned}
\end{equation}

We optimize the hyperparameters $\myVec{\theta}_{t}$ based on  \eqref{eqn:mseloss-sequntial} using deep learning optimization techniques based on, e.g., mini-batch stochastic gradient descent. 
In general, this form of learning based on first-order stochastic optimization requires one to be able to compute the gradient of \eqref{eqn:mseloss} and \eqref{eqn:mseloss-sequntial} with respect to the trainable parameters $\myVec{\theta}$ and $\myVec{\theta}_{t}$, respectively. By the formulation of \ac{dadmm}, and particularly, \eqref{eqn:PrimalStep}-\eqref{eqn:DualStep}, we note that the output of \ac{dadmm} with $T$ iterations is indeed differentiable with respect to these hyperparameters, as long as the local objectives $\{f_p(\cdot)\}$ are differentiable (which is implicitly assumed when formulating \ac{dadmm} with gradient steps, as in Algorithm~\ref{alg:Algo1}). 
The resulting procedure  is summarized as Algorithm~\ref{alg:L2O} assuming mini-batch stochastic gradient descent, which can naturally be extended to alternative learning mechanisms based on momentum and/or adaptive step-sizes~\cite{kingma2014adam}.
We initialize $\myVec{\theta}$  before the training process with fixed hyperparameters with which \ac{dadmm} converges (without a limited communications budget). 
After training,   the learned parameters  $\myVec{\theta}$  are used as hyperparameters for applying \ac{dadmm} with $T$ iterations.

  \begin{algorithm}
    \caption{Training Unfolded \ac{dadmm}}
    \label{alg:L2O} 
    \SetKwInOut{Initialization}{Init}
    \Initialization{Set $\myVec{\theta}$ as fixed hyperparameters \newline Fix learning rate $\mu>0$ and  epochs $i_{\max}$}
    \SetKwInOut{Input}{Input}  
    \Input{Training set  $\mathcal{D}$}  
    {
        \For{$i = 0, 1, \ldots, i_{\max}-1$}{%
                    Randomly divide  $\mathcal{D}$ into $Q$ batches $\{\mathcal{D}_q\}_{q=1}^Q$\;
                    \For{$q = 1, \ldots, Q$}{
                   Apply \ac{dadmm} with parameters $\myVec{\theta}$ to $\mathcal{D}_q$\;
                    
                    Compute batch loss $\mathcal{L}_{\mathcal{D}_q}(\myVec{\theta})$  by \eqref{eqn:mseloss}\;
                    Update  $\myVec{\theta}\leftarrow \myVec{\theta} - \mu\nabla_{\myVec{\theta}}\mathcal{L}_{\mathcal{D}_q}(\myVec{\theta})$\;
                    }
                    
                    }
        \KwRet{$\myVec{\theta}$}
  }
\end{algorithm}


\subsection{Discussion}
\label{ssec:discussion}
The proposed  unfolded \ac{dadmm} leverages data, obtained from simulations or measurements, to enable the \ac{dadmm} algorithm to operate reliably and rapidly with a fixed number of communication rounds. Distributed optimization algorithms are typically evaluated and analyzed in terms of their convergence, which is an asymptotic property. Here, we are particularly interested in operation with a fixed and typically small number of communication rounds (e.g., in our numerical studies we use $T = 20$), acknowledging that often in practice it is highly beneficial to impose such limitations on distributed systems. For this reason, we design unfolded \ac{dadmm} to optimize the \ac{dadmm} distributed optimizer within this predefined finite horizon.
\textcolor{NewColor}{The selection of $T$ depends on the communication latency that the application initiating distributed optimization is willing to tolerate, and thus varies between applications. Accordingly, we treat $T$ is a fixed parameter provided during design based on the application and system considerations.}

While we use deep learning techniques to optimize the optimizer, the resulting distributed inference rule is not a black-box \ac{dnn}, and maintains the interpretable operation of the classic \ac{dadmm} algorithm. This follows since deep learning automated training machinery is utilized to adapt the hyperparameters of the algorithm and the regularization coefficient of the objective, which are typically tuned by hand, while allowing these parameters to vary between iterations. 
\textcolor{NewColor}{Once trained, unfolded \ac{dadmm} implements exactly $T$ rounds of conventional \ac{dadmm} optimization, operating with the same complexity of its model-based counterpart applied with fixed hyperparameters.
The fact that the performance of iterative optimizers can be improved by proper setting of its hyperparameters is well established in the optimization literature~\cite{boyd2004convex}; our unfolded design} removes the need to adapt these hyperparameters manually or via additional lengthy computations during inference (via e.g., backtracking or line-search \cite[Ch. 9.2]{boyd2004convex}). This allows to notably improve both accuracy and run-time, as experimentally demonstrated in Subsections~\ref{ssec:lasso-results} and~\ref{ssec:linear-results} for the case studies of distributed sparse recovery  and of distributed linear regression learning, respectively. 

\textcolor{NewColor}{Our proposed unfolded \ac{dadmm} is designed for distributed optimization problems that one can tackle using conventional \ac{dadmm}. Such problems require  one to be able to compute the objective gradients in \eqref{eqn:PrimalStep}-\eqref{eqn:DualStep}, and assume that the optimized variable $\bar{\myVec{y}}$ takes values in a continuous set. When the objective is non-differentiable,  one potential solution is to seek a differentiable proxy that can approximate the non-differentiable objective function, and use this proxy during the training phase~\cite{shlezinger2022deep}. Similarly, when optimizing over discrete sets, one can train the unfolded optimizer while relaxing the need for the optimized variable to take discrete values, as in~\cite{khobahi2021lord}. As our work focuses on distributed optimization problems that can are naturally tackled using \ac{dadmm}, we leave the exploration of alternative training approaches  for future work.}

\textcolor{NewColor}{We propose two approaches to parameterize \ac{dadmm} when unfolded into a machine learning architecture in 
Subsection~\ref{ssec:Architecture}. The agent specific architecture provides increased flexibility, yet its number parameters grows with the number of agents, which in some settings may affect complexity and scalability. This is alleviated in the shared hyperparameters design, where the number of parameters is invariant of the number of agents. The  constrained form of shared hyperparameters design comes at the cost of a relatively minor performance degradation, while facilitating scalability and making training less complex by allowing an architecture trained on a small graph to be used for reliable distributed optimization over large graphs, as empirically observed in Subsection~\ref{ssec:linear-results}. One can consider alternative architectures that achieve different balances between parameterization, abstractness, and training complexity. For instance, a possible hybrid approach, which is left for future investigation, has clusters of agents sharing hyperparameters~\cite{shlezinger2020communication}. 
Alternatively, training complexity can be affected by using loss measures that account for intermediate iterations as in~\cite{samuel2019learning}, leveraging the fact that these features can be associated with gradually improved estimates. Moreover, the fact that estimates can be evaluated using the optimization objective in \eqref{eqn:AugLag} indicates that our approach can be utilized for unsupervised learning~\cite{lavi2023learn}. We leave these extensions to future work.}

The unfolding of a distributed optimizer unveils a key benefit of deep unfolding in its ability to reduce not only the run-time of iterative optimizers as in conventional centralized deep unfolding~\cite{shlezinger2022model}, but also lead to savings in communications. Our formulation considers the unfolding of a specific distributed optimization algorithm,  \ac{dadmm}~\cite{mota2013d}, \color{NewColor}which extends the centralized \acs{dadmm} optimizer to distributed settings~\cite{boyd2011distributed}. While various works combined centralized \acs{admm} with deep unfolding~\cite{johnston2021admm,zhou2021admm,wang2022efficient,fan2023fast,shah2024optimization}, our design, which is based on \ac{dadmm}, is notably different from unfolding of centralized optimizers. In addition to its parameterization of the distributed optimizer, our design also incorporates the need to account for the presence of multiple agents and underlying graphs that are not present in centralized settings (yet affect our different designs detailed in Subsection~\ref{ssec:Architecture});  most importantly, it reveals the core gains of deep unfolding for distributed optimization in terms of facilitating operation with limited communications.

We formulate the unfolded distributed optimizer for a generic distributed optimization setup, 
and specialize it for the representative case studies  distributed sparse recovery and linear regression. 
\color{black}
However,  the consistent improvements in performance and run-time of our unfolded \ac{dadmm} reported in the sequel indicates on the potential of this approach to other optimization problems and alternative distributed algorithms, which can possibly operate with fewer communications when combined with our considered form of deep unfolding. However, we leave the generalization of this study to other forms of distributed optimization for future work. 
%
Finally, while our derivation does not distinguish between different links, one can extend our approach to settings where some links are more constrained or costly by using weighted graphs, as well as be combined with distributed optimization with quantized messages. These extensions are all left for future study.

%
\section{Case Study 1: Distributed Sparse Recovery}
\label{sec:LASSO-Model}
%
%
The formulation of our proposed unfolded \ac{dadmm} in Section~\ref{sec:Unfolded} derives the algorithm for a generic distributed optimization setup. However, its conversion of \ac{dadmm} into a trainable machine learning model highly depends on the specific problem considered. Therefore, we next provide two case studies for which we specialize and evaluate unfolded \ac{dadmm}: A distributed sparse recovery, detailed in this section, and distributed linear regression  presented in Section~\ref{sec:LinearRegModel}. 

We particularly focus on the application of unfolded \ac{dadmm} for distributed sparse recovery based on the  \ac{dlasso} objective, which is a  convex relaxation based problem formulation commonly employed in compressed sensing setups~\cite[Ch. 1]{eldar2012compressed}. In Subsection~\ref{ssec:lasso-problem} we describe the \ac{dlasso} problem formulation. Then, we specialize unfolded \ac{dadmm} for solving this sparse recovery problem in Subsection~\ref{ssec:sparse-recovery} and experimentally evaluate it in Subsection~\ref{ssec:lasso-results}.


\subsection{D-LASSO Problem Formulation}
\label{ssec:lasso-problem} 
%

We consider the recovery of a sparse signal from a set of compressed measurements observed by the agents. 
Accordingly, the local measurements represent compressed versions of some sparse high dimensional vector obtained using a set of sensing matrices. The distributed nature of the problem arises from the decentralized data residing on the computing agents. 

The \ac{dlasso} problem formulates this task by relaxing it into a $\ell_1$ regularized recovery objective.  Here, the functions $\{f_p(\cdot)\}$ in \eqref{eqn:opt-problem} are given by
\begin{equation}
    \label{eqn:LassoObj}
    f_p(\myVec{\Bar{y}}; \myVec{b}_p) = \frac{1}{2}\|\myMat{A}_p\myVec{\Bar{y}} - \myVec{b}_p\|_2^2 + \tau \|\myVec{\Bar{y}}\|_1,
\end{equation}
for each $p \in \mySet{V}$. In \eqref{eqn:LassoObj}, $\myMat{A}_p\in\mathbb{R}^{m\times n}$ is the sensing matrix with which the $p$th agent acquires its local observation $\myVec{b}_{p} \in \mathbb{R}^{m}$,  with $m < n$. The regularization coefficient $\tau>0$ balances the sparsity of the solution and its matching of the data, and is a hyperparameter originating from the representation of the sparse recovery task as a relaxed $\ell_1$ regularized objective. This objective hyperparameter is often tuned manually.


%
\subsection{Unfolded D-ADMM for the \ac{dlasso}}
\label{ssec:sparse-recovery} 
%

\ac{dadmm} can be utilized for tackling the D-LASSO problem in \eqref{eqn:LassoObj}, which we in turn unfold following the methodology detailed in Section~\ref{sec:Unfolded}. Therefore, in the following we first specialize \ac{dadmm} for the \ac{dlasso} problem, after which we formulate its associated unfolded \ac{dadmm} machine learning model. 

\subsubsection{\ac{dadmm} for the \ac{dlasso}}
As described in Subsection~\ref{ssec:d-admm}, \ac{dadmm} converts the objective via variable splitting into
\begin{align}
 \label{eqn:distopt-problem}
        &\mathop{\arg \min}\limits_{\myVec{y}_{1}, \ldots, \myVec{y}_{P}
    } \quad  \sum_{p=1}^{P} \frac{1}{2}\|\myMat{A}_p\myVec{{y}}_p - \myVec{b}_p\|_2^2 + \tau \|\myVec{{y}}_p\|_1, \\
    &\text{subject to}  \quad \myVec{y}_{p} = \myVec{y}_{j}, \forall j \in \mySet{N}_p,\notag
    \end{align}
    where $\mySet{N}_p$ is the set of neighbors of node $p$ in  $G(\mySet{V}, \mySet{E})$.
    The augmented Lagrangian based on \eqref{eqn:distopt-problem}  for the $p$th agent is formulated as
\begin{align}
    &\mathcal{L}_p(\myVec{y}_{p}, \{\myVec{y}_{j}\}_{j\in\mySet{N}_p}, \myVec{\lambda}_p) \triangleq 
    \frac{1}{2}\|\myMat{A}_p\myVec{{y}}_p - \myVec{b}_p\|_2^2 + \tau \|\myVec{{y}}_p\|_1 \notag \\
    &\qquad \qquad + \sum_{j \in \mySet{N}_p} \myVec{\lambda}_p^T(\myVec{y}_{p} - \myVec{y}_{j}) + \frac{\rho}{2}\|\myVec{y}_{p} - \myVec{y}_{j}\|_2^2,
    \label{eqn:LASSOAugLag}
\end{align}
where $\rho>0$ is a fixed hyperparameter and $\myVec{\lambda}_p$ is the dual variable of node $p$. 

Accordingly, \ac{dadmm} (Algorithm~\ref{alg:Algo1}) can be applied to tackling the \ac{dlasso} problem in \eqref{eqn:LassoObj} with \eqref{eqn:PrimalStep} becoming
\begin{align}
    \myVec{y}_{p}^{(k+1)} &= \myVec{y}_{p}^{(k)} -\alpha \nabla_{\myVec{y}_{p}} \mathcal{L}_p(\myVec{y}_{p}^{(k)}, \{\myVec{y}_{j,p}\}_{j\in\mySet{N}_p}, \myVec{\lambda}_p^{(k)}) \notag \\
    &=\myVec{y}_{p}^{(k)} -\alpha \Big( \myMat{A}_p^T\myMat{A}_p\myVec{{y}}_p^{(k)} - \myMat{A}_p^T\myVec{b}_p^{(k)} + \tau \cdot {\rm sign}\big(\myVec{{y}}_p^{(k)}\big) \notag \\
    & \qquad\qquad  +  \sum_{j \in \mySet{N}_p} \myVec{\lambda}_p^{(k)} + \rho\big(\myVec{y}_{p}^{(k)} - \myVec{y}_{j,p} \big) 
    \Big),
    \label{eqn:LASSOPrimalStep}
\end{align}
by each agent $p$, with $\alpha>0$ being a step-size.
Similarly, \eqref{eqn:DualStep} is specialized into
\begin{align}
    \myVec{\lambda}_{p}^{(k+1)} &= \myVec{\lambda}_{p}^{(k)} +\eta \nabla_{\myVec{\lambda}_{p}} \mathcal{L}_p(\myVec{y}_{p}^{(k+1)}, \{\myVec{y}_{j,p}\}_{j\in\mySet{N}_p}, \myVec{\lambda}_p^{(k)}) \notag \\
    &= \myVec{\lambda}_{p}^{(k)} +\eta\sum_{j \in \mySet{N}_p}\big(\myVec{y}_{p}^{(k+1)} - \myVec{y}_{j,p}\big),
    \label{eqn:LASSODualStep}
\end{align}
 with step-size $\eta>0$ being the dual variable step size.

\begin{figure}
    \centering
    \includegraphics[width=\columnwidth]{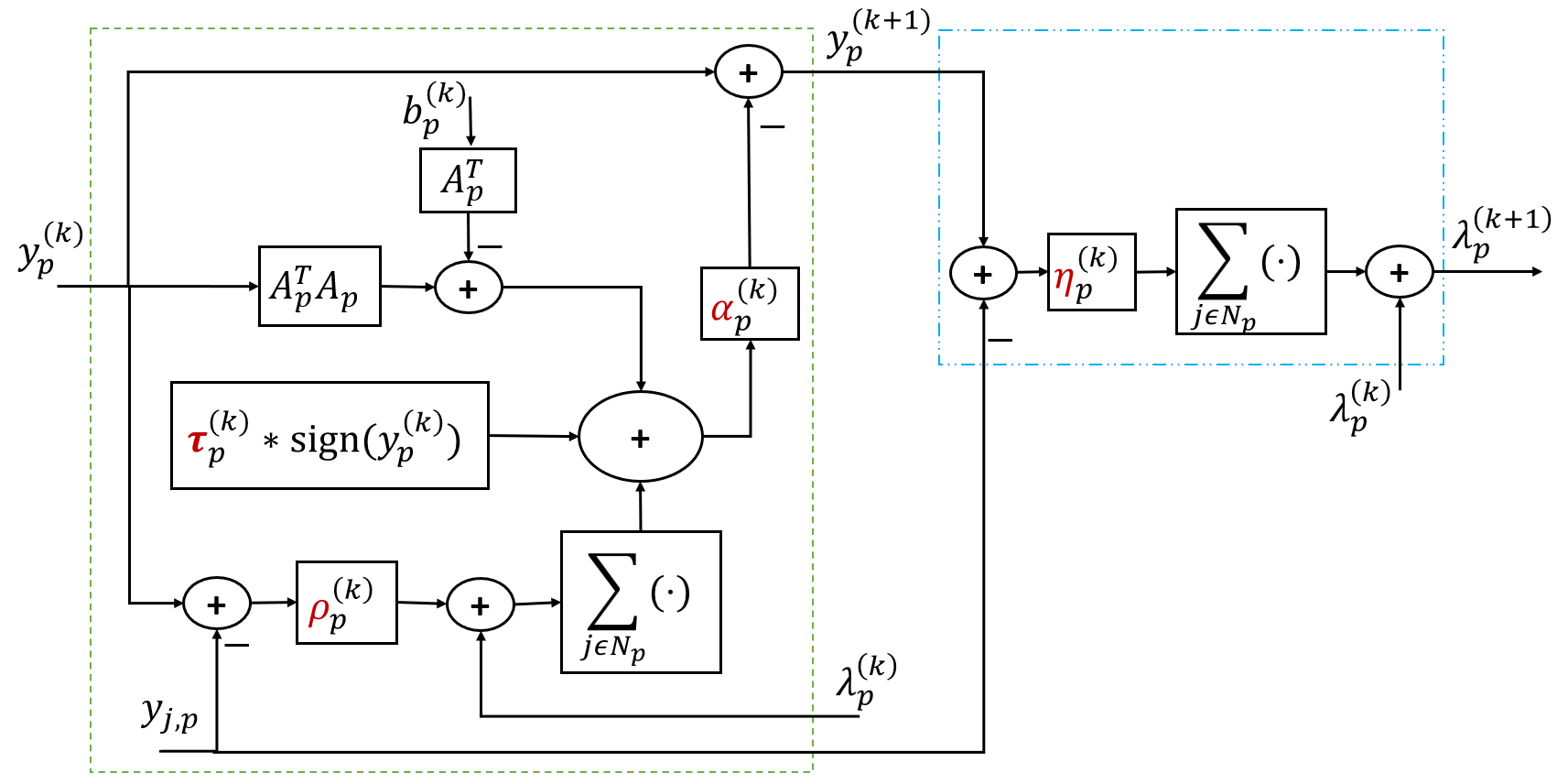}
    \caption{Unfolded \ac{dadmm} for LASSO at agent $p$ in iteration $k$. Dashed green and blue blocks are the primal and dual updates, respectively. Red fonts represent trainable parameters.}
    \label{fig:primal-block-diagram}
\end{figure}

\subsubsection{Unfolding \ac{dadmm}}
Following the description in Section~\ref{sec:Unfolded}, we unfold \ac{dadmm} into a machine learning model by fixing its number of iterations to be $T$ and treating the hyperparameters of \ac{dadmm}, i.e., $ \alpha, \rho,$ and $\eta$ as trainable parameters. Moreover, we note that the formulation of the \ac{dlasso} objective introduces an additional hyperparameter, $\tau$, which is not unique to its tackling by \ac{dadmm}. Since the tuning of this parameter can largely affect the performance of iterative optimizers applied to such objectives, we leverage our available data to also treat $\tau$ as a trainable parameter.

Accordingly, the trainable parameters of unfolded \ac{dadmm} are given by the set of iteration-specific hyperparameters, i.e., $\myVec{\theta}=\{\myVec{\theta}_k\}_{k=1}^T$. The resulting trainable architecture is illustrated in Fig.~\ref{fig:primal-block-diagram}. When applying an agent-specific hyperparameters, the trainable parameters are  $\myVec{\theta}_{k} = \{\rho_{p}^{(k)}, \alpha_{p}^{(k)}, \eta_{p}^{(k)}, \tau_{p}^{(k)}\}_{p=1}^{P}$, i.e.,  $4\cdot P \cdot T$  parameters.

\subsection{Numerical Evaluation}
\label{ssec:lasso-results} 
%


We numerically evaluate the proposed  unfolded D-ADMM algorithm\footnote{The source code and hyperparameters used in this experimental study as well as the one reported in Subsection~\ref{ssec:lasso-results} are available  at \url{https://github.com/yoav1131/Deep-Unfolded-D-ADMM.git}.}, comparing it to the D-ADMM algorithm~\cite{mota2013d}, where we used fixed hyperparameters manually tuned based on empirical trials to systematically achieve convergence. We also compare unfolded \ac{dadmm} with a data-driven \ac{gnn}, trained with the same dataset, where we employ a \ac{gnn} based on the popular GraphSage architecture~\cite{hamilton2017inductive}. 

\begin{figure}
    \centering
    \includegraphics[width=\columnwidth]{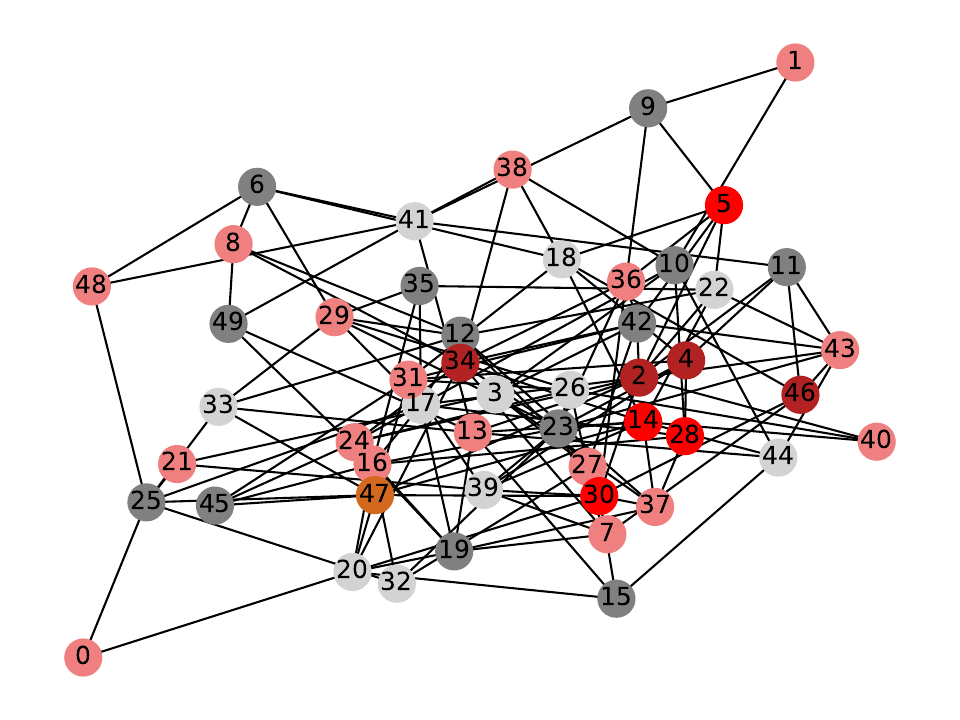}
    \caption{Proper colored network with $P=50$ nodes example.} 
    \label{fig:ProperColoredGraphP=50}
\end{figure}

We simulate a communication network using the Erdos-Renyi graph model of $P \in \{5, 20, 50\}$ nodes with proper coloring. An example of a graph generated with $P=50$ vertices is illustrated in Fig.~\ref{fig:ProperColoredGraphP=50}.
We generate observations for each node as $\myVec{b}_{p} = \myMat{A}_{p}\myVec{\bar{y}} + \myVec{n}_{p}$ with $n=2000$ and $m=\frac{500}{P}$, where $\myVec{n}_{p}$ is Gaussian  with zero-mean i.i.d. entries of variance $\sigma^2$. The desired $\myVec{\bar{y}}$ has $25\%$ non-zero entries, and the sensing matrices $\{\myMat{A}_{p}\}$ are taken to be the block sub-matrices of the full sensing matrix is taken from Problem 902 of the Sparco toolbox~\cite{friedlander2012dual}. 
The \ac{snr}, defined as $1/\sigma^2$, is set to  ${\rm SNR} \in \{-2, 0, 2, 4\}$ [dB].

\begin{figure*}
    \centering
    \begin{subfigure}{0.32\textwidth}
    \includegraphics[width=\linewidth]{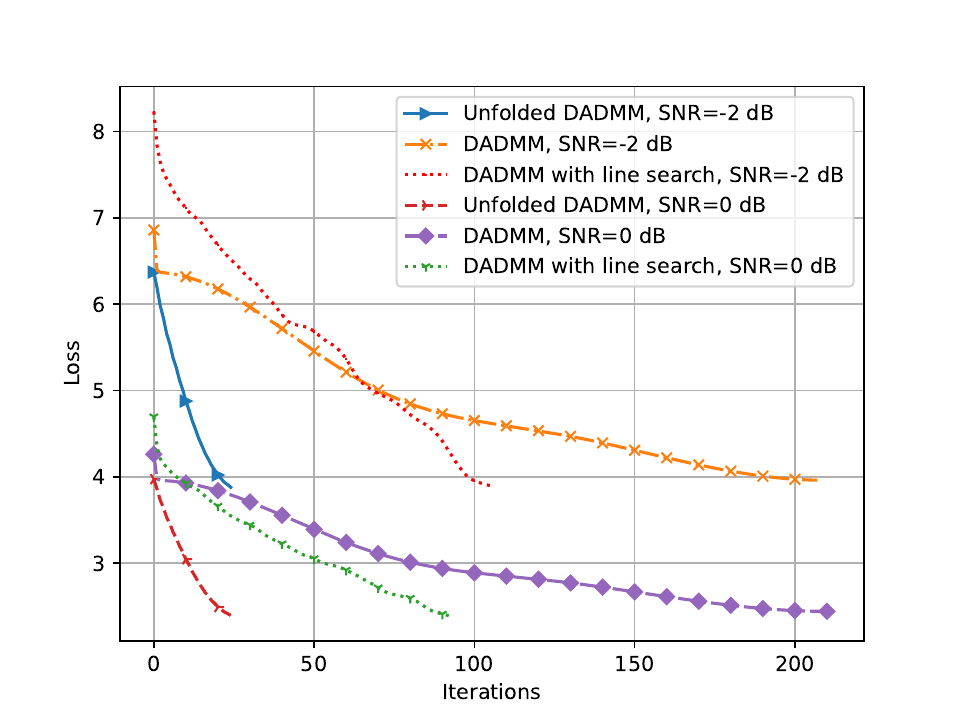}
     \caption{${\rm SNR} =-2,0$ dB, $P=5$}
     \label{fig:loss-P=5-1}
    \end{subfigure}
    \begin{subfigure}{0.32\textwidth}
    \includegraphics[width=\linewidth]{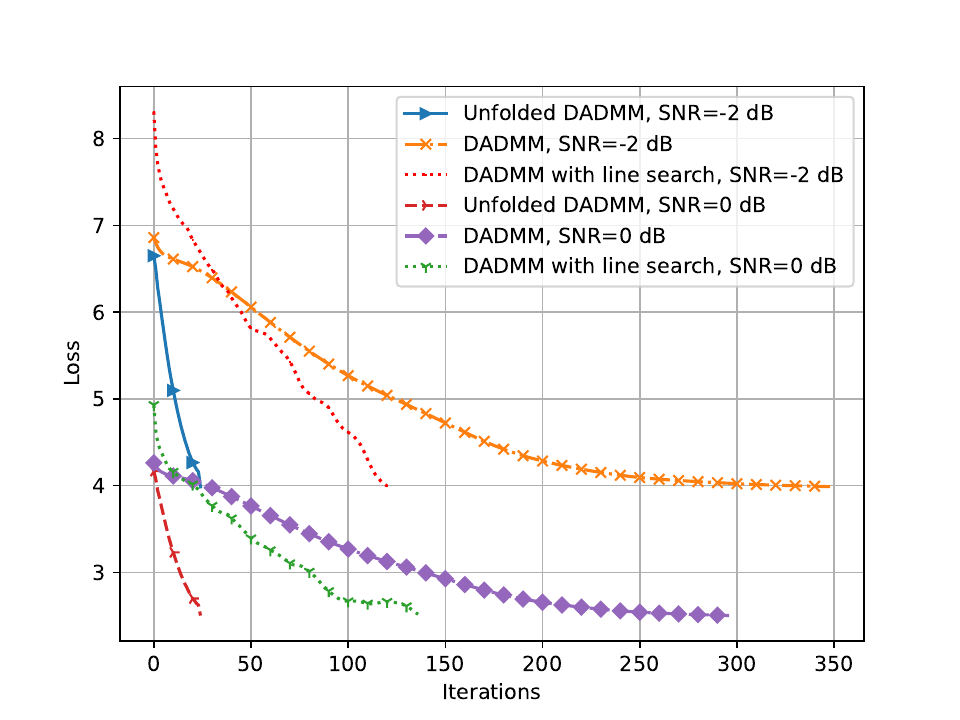}
     \caption{${\rm SNR} =-2,0$ dB, $P=20$}
     \label{fig:loss-P=20-1}
    \end{subfigure}
    \begin{subfigure}{0.32\textwidth}
    \includegraphics[width=\linewidth]{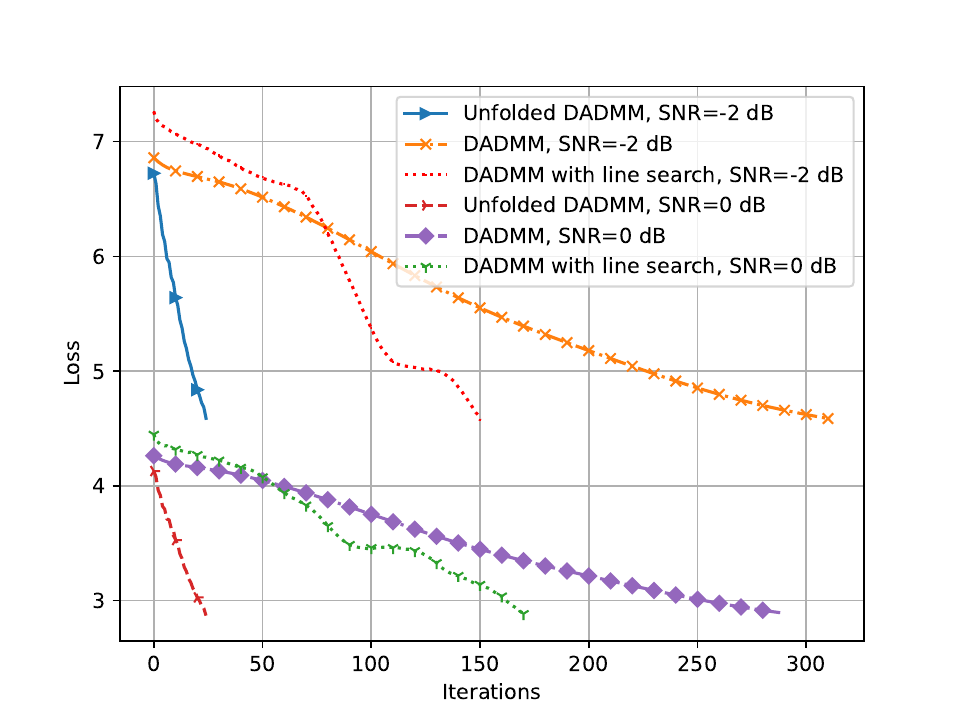}
     \caption{${\rm SNR} =-2,0$ dB, $P=50$}
     \label{fig:loss-P=50-1}
     \end{subfigure}
    \\
    \begin{subfigure}{0.32\textwidth}
    \includegraphics[width=\linewidth]{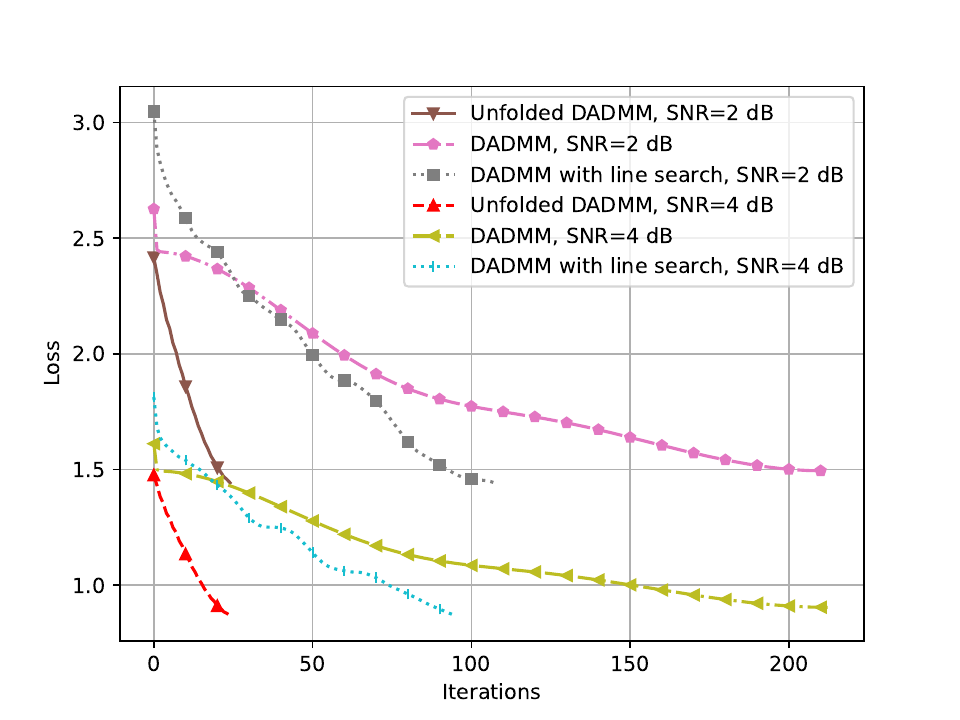}
    \caption{${\rm SNR} =2,4$ dB, $P=5$}
     \label{fig:loss-P=5-2}
    \end{subfigure}
    \begin{subfigure}{0.32\textwidth}
    \includegraphics[width=\linewidth]{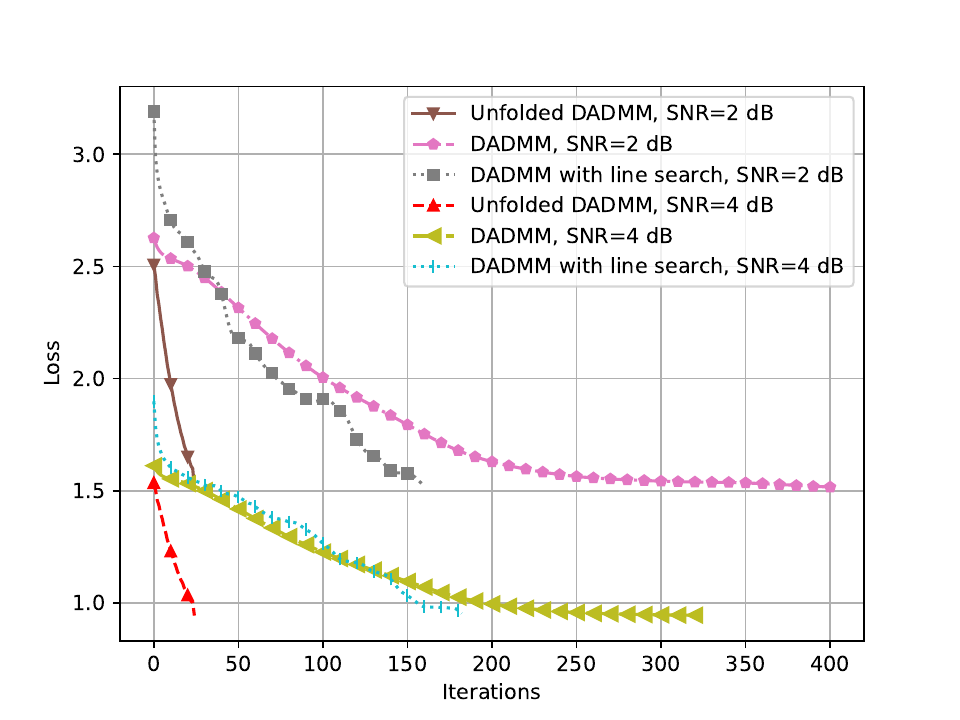}
    \caption{${\rm SNR} =2,4$ dB, $P=20$}
    \label{fig:loss-P=20-2}
    \end{subfigure}
    \begin{subfigure}{0.32\textwidth}
    \includegraphics[width=\linewidth]{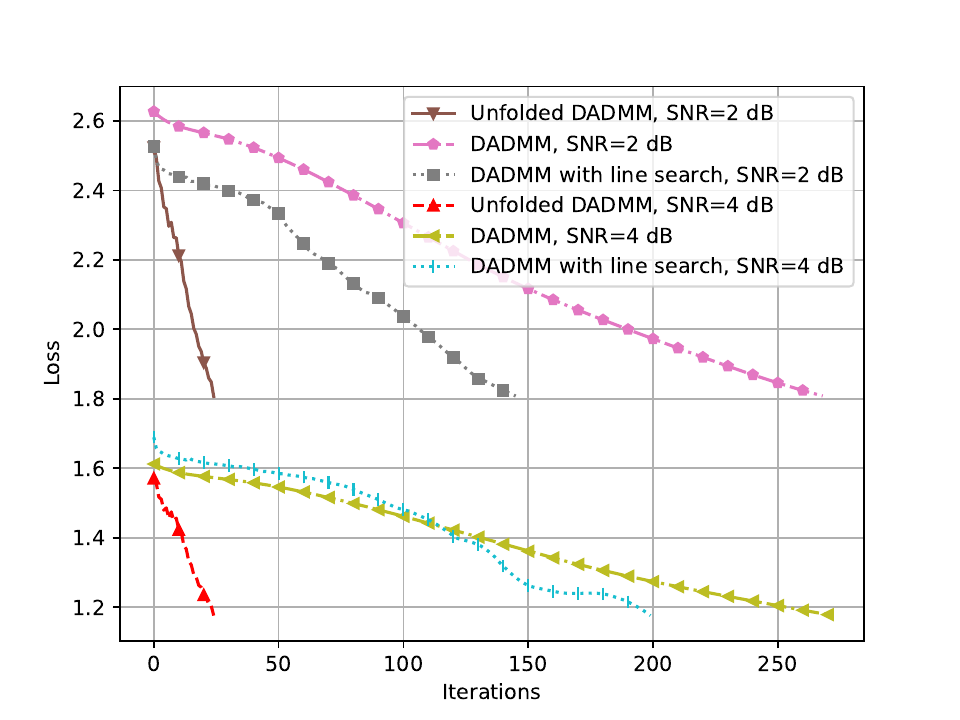}
    \caption{${\rm SNR} =2,4$ dB, $P=50$}
     \label{fig:loss-P=50-2}
     \end{subfigure}
    \caption{Loss per iteration, distributed LASSO.}
    \label{fig:lossVsIter}
\end{figure*}

For each simulated setting, we apply Algorithm~\ref{alg:L2O} to optimize $\myVec{\theta}$ with $T=25$ iterations based on a labeled dataset as in \eqref{eqn:dataset} comprised of $L=1200$ training samples. Both unfolded \ac{dadmm} and the \ac{gnn} are trained  using $100$ epochs of Adam \cite{kingma2014adam} with a batch size of $100$. 
All considered  algorithms are evaluated over $200$ test observations. 





\begin{figure*}
    \centering
    \begin{subfigure}{0.32\textwidth}
    \includegraphics[width=\linewidth]{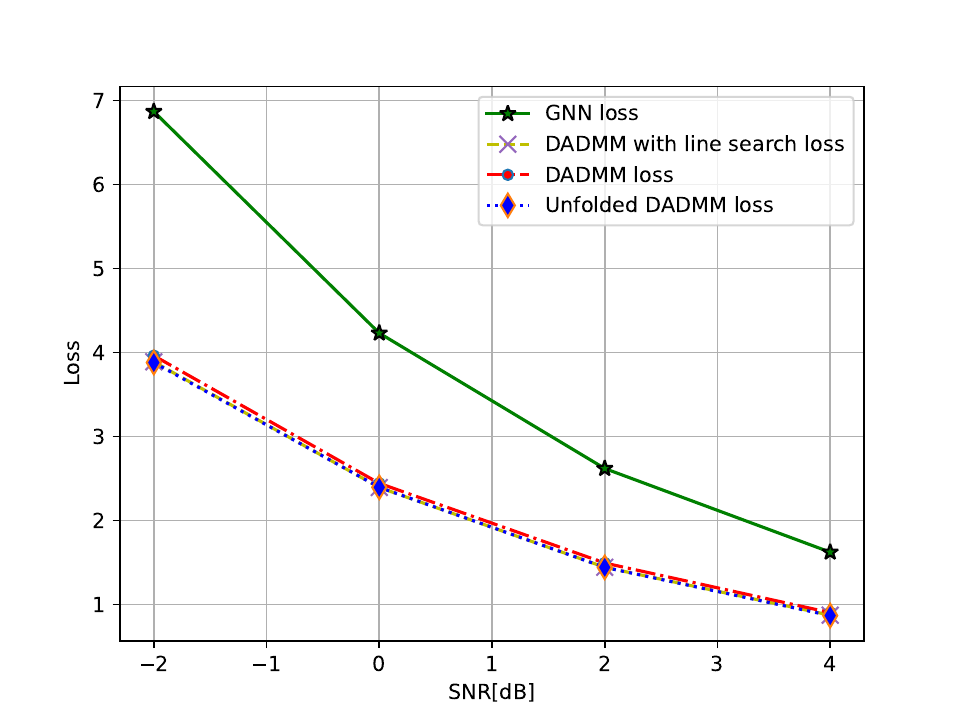}
    \caption{$P=5$}
    \label{fig:lossVSsnr-P=5}
    \end{subfigure}
    \begin{subfigure}{0.32\textwidth}
    \includegraphics[width=\linewidth]{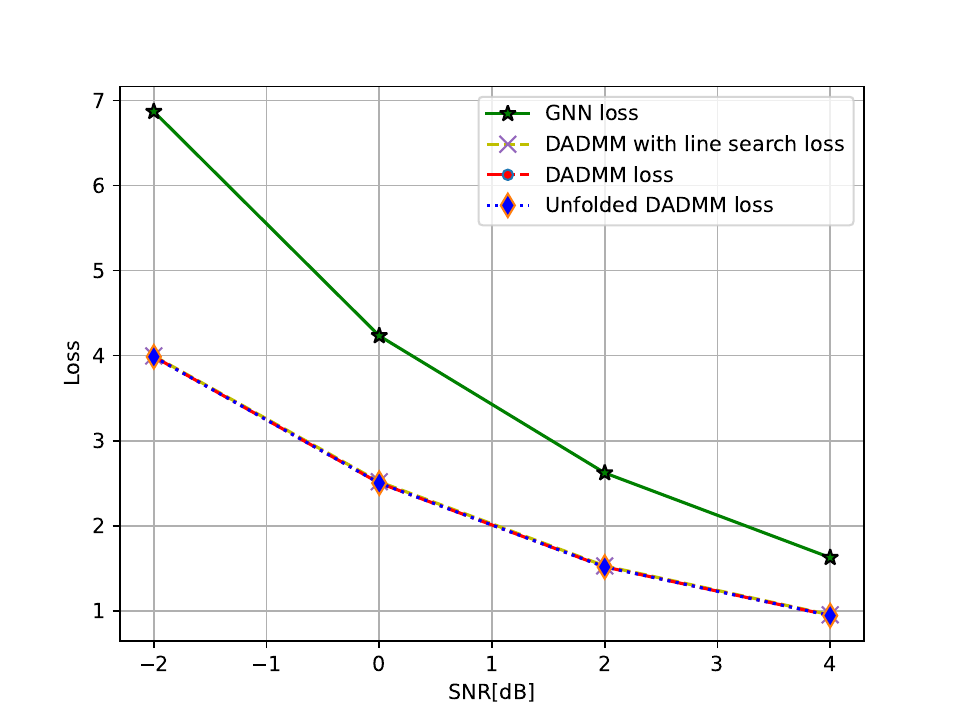}
    \caption{$P=20$}
    \label{fig:lossVSsnr-P=20}
    \end{subfigure}
    \begin{subfigure}{0.32\textwidth}
    \includegraphics[width=\linewidth]{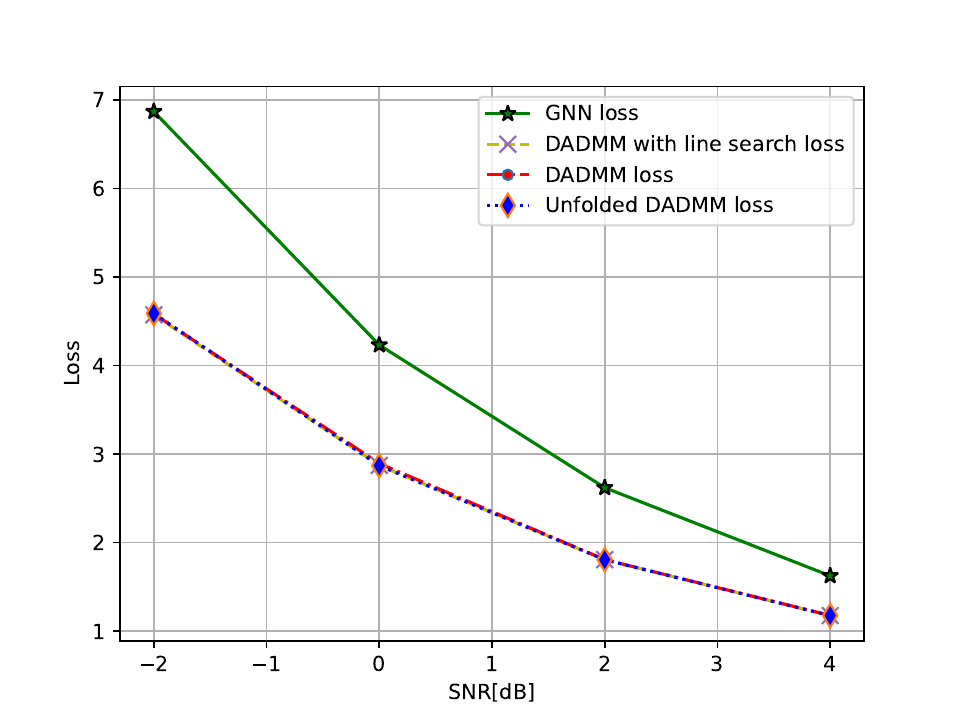}
    \caption{$P=50$}
    \label{fig:lossVSsnr-P=50}
    \end{subfigure}
    \caption{Final loss versus \ac{snr}, distributed LASSO.}
    \label{fig:LossVSsnr}
\end{figure*}

\begin{figure*}
    \centering
    \begin{subfigure}{0.32\textwidth}
    \includegraphics[width=\linewidth]{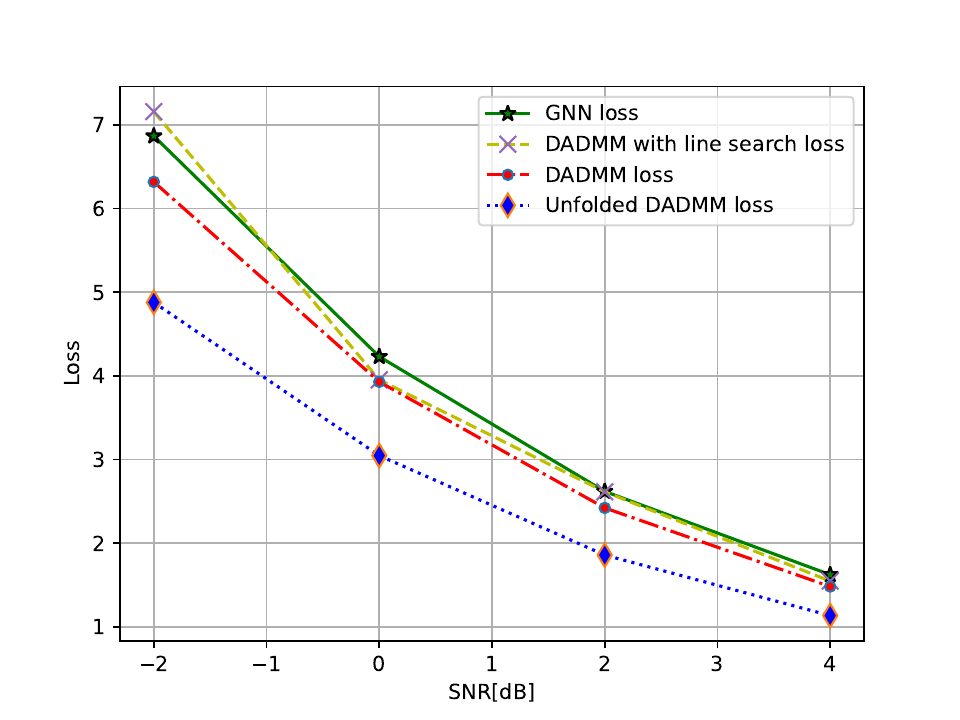}
    \caption{$P=5$}
    \label{fig:lossVSsnr10-P=5}
    \end{subfigure}
    \begin{subfigure}{0.32\textwidth} 
    \includegraphics[width=\linewidth]{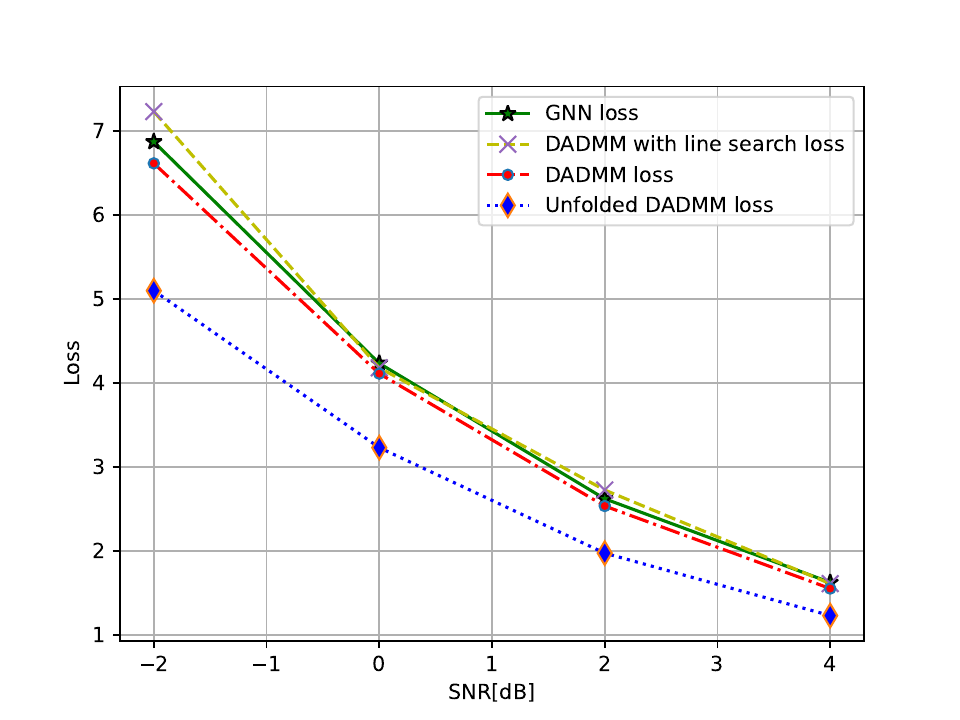}
    \caption{$P=20$}
    \label{fig:lossVSsnr10-P=20}
    \end{subfigure}
    \begin{subfigure}{0.32\textwidth}
    \includegraphics[width=\linewidth]{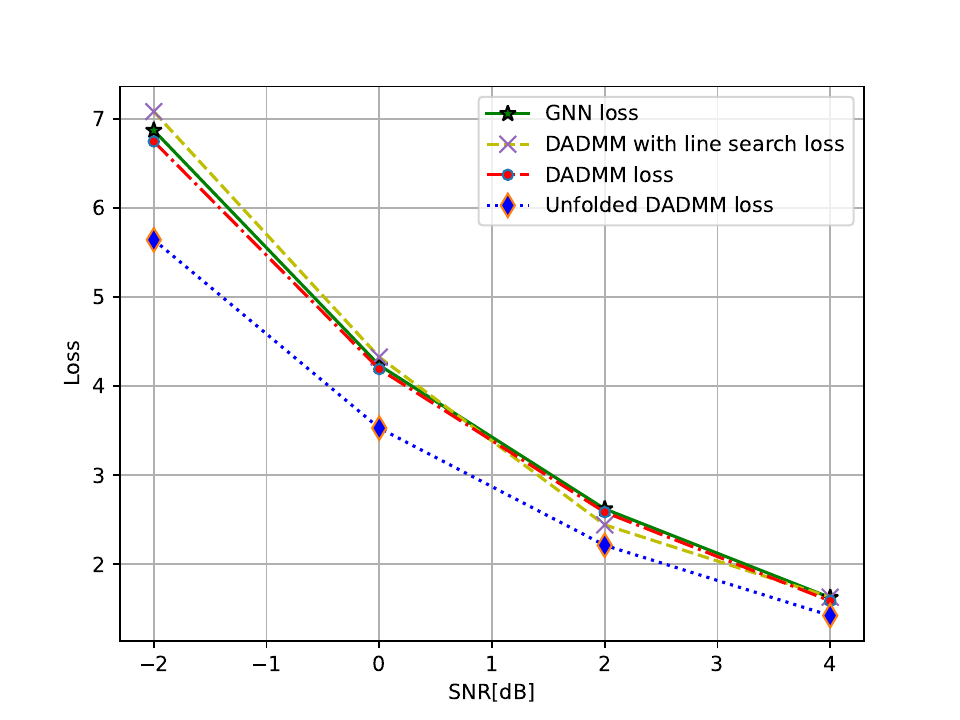}
    \caption{$P=50$}
    \label{fig:lossVSsnr10-P=50}
    \end{subfigure}
    \caption{Loss versus \ac{snr} after 10 communication rounds, distributed LASSO.}
    \label{fig:LossVSsnr10}
\end{figure*}



We first evaluate the loss versus iteration achieved by unfolded \ac{dadmm}  with agent specific hyperparameters compared with conventional \ac{dadmm}. The results are depicted in Fig.~\ref{fig:lossVsIter}. It is observed in Fig.~\ref{fig:lossVsIter} that the proposed unfolded D-ADMM improves upon the D-ADMM with fixed hyperparameters and \ac{dadmm} with line search optimization, by requiring much fewer iterations (hence communications) to converge. In particular,  for $P=5$ we observe in Figs.~\ref{fig:loss-P=5-1} and \ref{fig:loss-P=5-2} an average of communications reduction by factors of $\times 8$  and $\times 4$ ($25$ iterations vs. $211$ iterations and $101$ iterations, respectively), depending on the \ac{snr}. The corresponding reduction for $P=20$ is by $\times 13$ and $\times 6$ ($25$ iterations vs. $341$ iterations in Fig.~\ref{fig:loss-P=20-1} and $151$ iterations in Fig.~\ref{fig:loss-P=20-2}) while for $P=50$ we observe a reduction by $\times 11$ and $\times 6$ ($25$ iterations vs. $285$  iterations in Fig.~\ref{fig:loss-P=50-1} and $168$ iterations in Fig.~\ref{fig:loss-P=50-2}).

Next, we compare the performance achieved by unfolded \ac{dadmm} after its last iteration ($T=25$), compared with the \ac{gnn} with $10$ layers trained for the same task, as well as the conventional \ac{dadmm} when allowed to run until convergence. We observe in Fig.~\ref{fig:LossVSsnr} that for all considered graph sizes, our unfolded \ac{dadmm} effectively coincides with the converged \ac{dadmm}, while using only $T=25$ iterations. The \ac{gnn}, which is invariant of the optimizer operation and just learns an abstract message passing as layers of a \ac{dnn}, is  outperformed by our unfolded \ac{dadmm}.
Since the \ac{gnn} only employs $10$ communication rounds (being comprised of $10$ layers), in Fig.~\ref{fig:LossVSsnr10} we compare the performance of all optimizers after $T=10$ communication rounds. There we systematically observe the improved performance achieved by the unfolded optimizer. 
%
%
%
These results demonstrate the benefits of the proposed approach in leveraging data to improve both performance and convergence speed while preserving the interpretability and suitability of conventional iterative optimizers.

\begin{table}
    \centering
    \resizebox{0.49\textwidth}{!}{
    \begin{tabular}{|l|*{3}{c|}}\hline
    \backslashbox{Method}{No. of agents (P)} & \makebox{\textbf{$P = 5$}} & \makebox{\textbf{$P = 20$}} & \makebox{\textbf{$P = 50$}}  \\ [3pt] 
    \hline 
    GNN & $2.42 [s]$ & $62.24 [s]$ & $112.10 [s]$ \\ [3pt]
    \ac{dadmm} & $2.30 [s]$ & $53.18 [s]$ & $87.40[s]$ \\ [3pt]
    Unfolded \ac{dadmm} & $2.30 [s]$ & $53.18 [s]$ & $87.40[s]$ \\ [3pt]
    \ac{dadmm} with line search & $4.69 [s]$ & $97.36 [s]$ & $303.69[s]$ \\ [3pt]
    \hline
  \end{tabular}}
  \caption{Average run-time results for $25$ iterations of distributed optimization with different number of agents ($P$).} 
  \label{tab:runtime}
  
\end{table}

\color{NewColor}{We conclude our evaluation by comparing the run-time performance achieved by unfolded \ac{dadmm} with $25$ iterations, compared with the \ac{gnn} with $25$ layers trained for the same task, as well as the conventional \ac{dadmm} with $25$ iterations with fixed hyperparmaeters  and with line search. 
To that aim, we compute the average run-time of distributed optimization carried out over $100$ Erdős-Rényi random graphs generated for each number of agents  $P=\{5,20,50\}$. All algorithms are implemented on the same platform, that uses   NVIDIA GeForce RTX 3080 GPU, where all links are treated identically. 
The average run-time, reported in Table~\ref{tab:runtime}, show that for all considered graph sizes, the latency of our unfolded \ac{dadmm} coincide with that of conventional \ac{dadmm}, as both implement the same steps. Both are faster than optimizing using a \ac{gnn} or when employing line search, with gain becoming more substantial as the graph size increases.} 
\color{black}


%
\section{Case Study 2: Distributed Linear Regression}
\label{sec:LinearRegModel}
%
%
Our second case study considers a distributed machine learning task. Collaborative learning of machine learning models is becoming increasingly popular, particularly due to its ability to alleviate privacy considerations associated with data sharing~\cite{gafni2022federated}. A common framework for distributed machine learning is \ac{fl}, which typically deals with distributed learning with centralized orchestration, i.e., where there is a central server with which all agents communicate directly and that  enforces consensus on each communication round. Nonetheless, recent explorations have also considered its extensions to purely decentralized networks~\cite{ayache2023walk}.  

We particularly focus on learning a linear model based on the \ac{dlr} objective, as formulated in Subsection~\ref{ssec:linearReg-problem}. The application of unfolded \ac{dadmm} for this learning task is then presented in Subsection~\ref{ssec:regression} and numerically evaluated in Subsection~\ref{ssec:linear-results}.


\subsection{D-LR Problem Formulation}
\label{ssec:linearReg-problem} 
%
 We consider the learning of a linear model from a dataset that is partitioned into subsets and distributed across multiple computing agents. The local data here represents labeled sets used for learning purposes. Accordingly, the \ac{dlr} problem specializes the generic formulation in \eqref{eqn:opt-problem}  by setting the objectives  $\{f_p(\cdot)\}$ to be 
\begin{equation}
    \label{eqn:LinearRegObj}
    f_p(\myVec{\Bar{y}}; \myVec{b}_{p}) = \frac{1}{2 \cdot  L_{p} }\sum_{(\myVec{x}_{i,p}, s_{i,p} )\in \myVec{b}_{p}} (\myVec{a}^{T}\myVec{x}_{i, p} + \omega - s_{i, p})^{2}.
\end{equation}
In \eqref{eqn:LinearRegObj}, the local data $\myVec{b}_p$ for agent  $p \in \mySet{V}$  is the set of $L_p$ labeled pairs written as $\myVec{b}_{p} = \{\myVec{x}_{i, p}, s_{i, p}\}_{i=1}^{L_p}$, where $s_{i,p}$ is the scalar label associated with the $d\times 1$ input vector $\myVec{x}_{i,p}$. The optimization variable $\bar{\myVec{y}}$ is a linear regression model written as $\myVec{\bar{y}} = \{\myVec{\bar{a}}, \bar{\omega}\}$ representing an affine transformation with parameters $\myVec{\bar{a}} \in \mathbb{R}^{d}$ and $\bar{\omega} \in \mathbb{R}$ denoting the regression vector and the bias, respectively. 


%
\subsection{Unfolded D-ADMM for the \ac{dlr}}
\label{ssec:regression} 
%



\ac{dadmm} can be utilized for tackling the \ac{dlr} problem in \eqref{eqn:LinearRegObj}, which we in turn unfold following the methodology detailed in Section~\ref{sec:Unfolded}. Therefore, in the following we first specialize \ac{dadmm} for the \ac{dlr} problem, after which we formulate the unfolded \ac{dadmm} machine learning model.

\begin{figure*}
    \centering
    \begin{subfigure}{0.48\textwidth}
    \includegraphics[width=\textwidth]{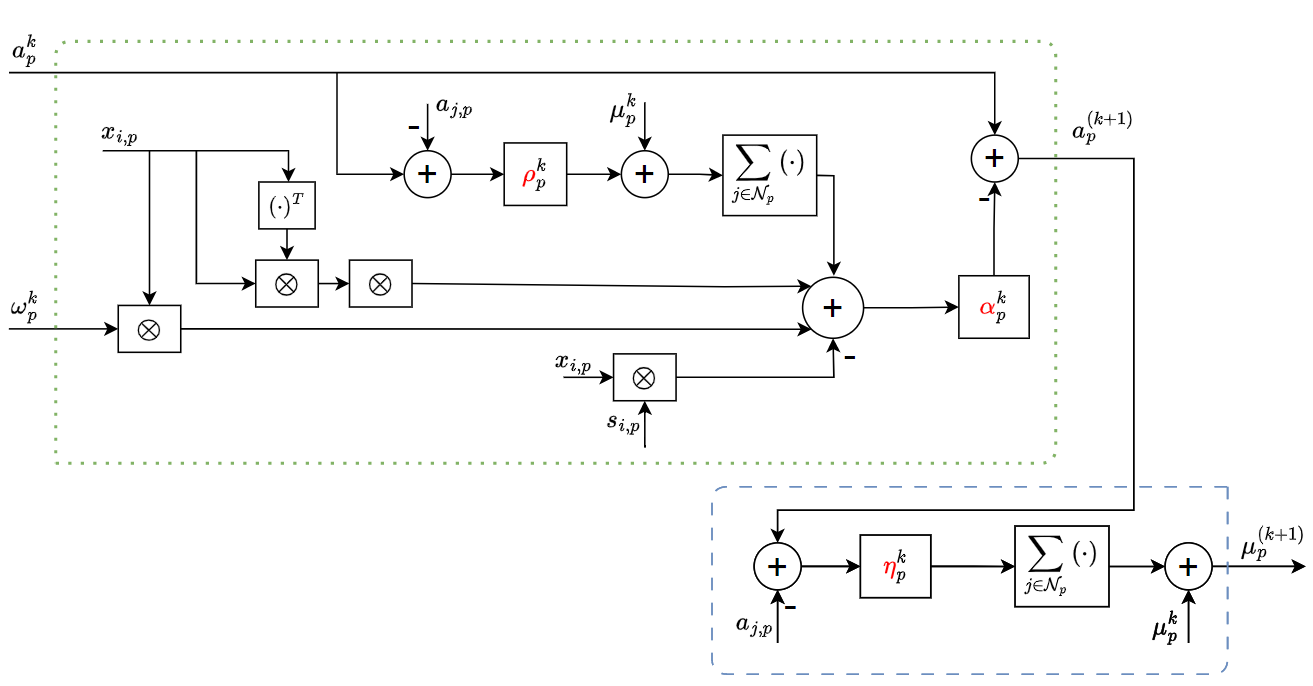}
    \caption{Optimized $\myVec{a}$ update.}
    \label{fig:LinearReg-block-diagram-a}
    \end{subfigure}
    \begin{subfigure}{0.48\textwidth}
    \includegraphics[width=\textwidth]
    {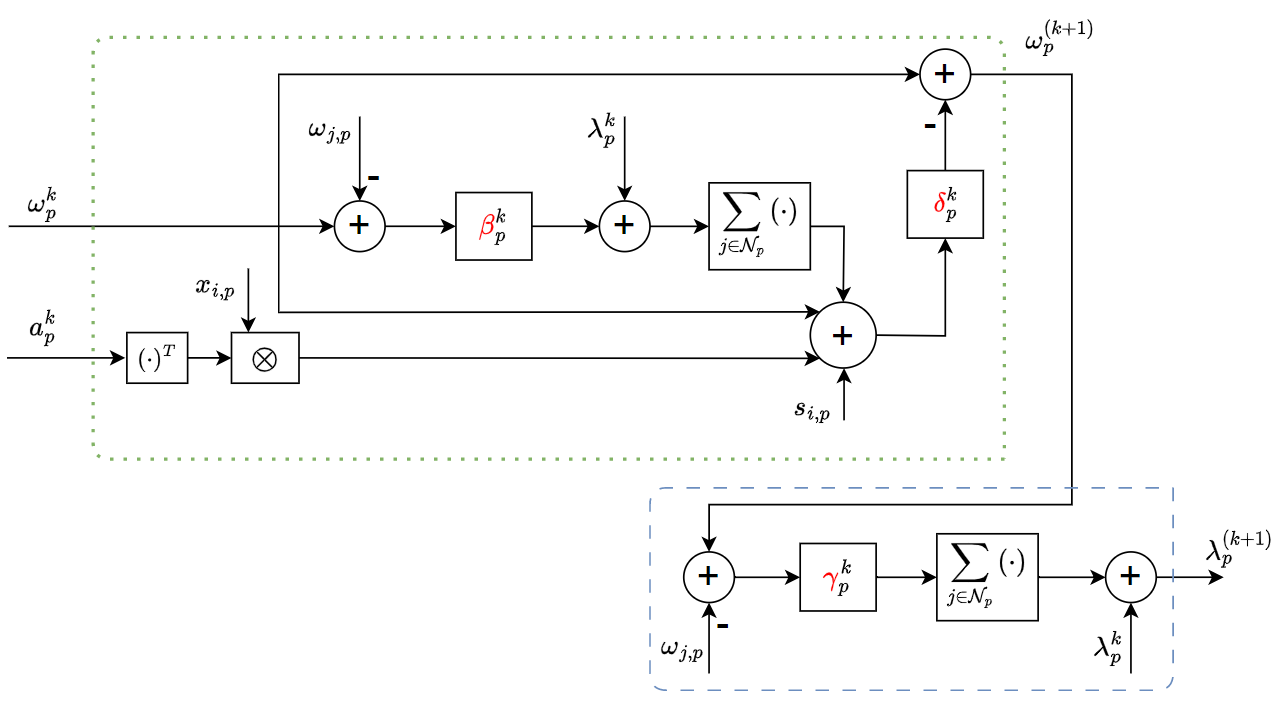}
    \caption{Optimized $\omega$ update.}
    \label{fig:LinearReg-block-diagram-omega}
    \end{subfigure}
    \caption{Unfolded \ac{dadmm} for linear regression model illustration at agent $p$ in iteration $k$. Dashed green and blue blocks are the primal update and the dual update, respectively. Red fonts represent trainable parameters.}
    \label{fig:LinearReg-block-diagram}
\end{figure*}

\subsubsection{\ac{dadmm} for the \ac{dlr}}
Following Subsection~\ref{ssec:d-admm}, \ac{dadmm} converts the objective via variable splitting into
\begin{align}
 \label{eqn:distlinear-problem}
        &\mathop{\arg \min}\limits_{\myVec{y}_{1}, \ldots, \myVec{y}_{P}
    } \quad \sum_{p=1}^{P}  \sum_{(\myVec{x}_{i,p}, s_{i,p}) \in \myVec{b}_{p}} \frac{1}{2  \cdot L_{p}}(\myVec{a}_{p}^{T}\myVec{x}_{i, p} + \omega_{p} - s_{i,p})^{2}, \\
    &\text{subject to}  \quad \myVec{a}_{p} = \myVec{a}_{j},  \forall j \in \mySet{N}_p,\notag \\
    &\qquad \qquad \quad \omega_{p} = \omega_{j}, \forall j \in \mySet{N}_p,\notag
    \end{align}
    where the local models 
    $\myVec{y}_{p}$ have two components, $\{\myVec{a}_{p}, \omega_{p}\}$, hence \eqref{eqn:distlinear-problem} has two constraints. This leads to two primal and two dual problems one for each component. 
    The augmented Lagrangian based on \eqref{eqn:distlinear-problem}  for the $p$th agent is formulated as
\begin{align}
    &\mathcal{L}_p(\myVec{y}_{p}, \{\myVec{y}_{j}\}_{j\in\mySet{N}_p}, \myVec{\mu}_{p}, \myVec{\lambda}_{p}) \triangleq 
    \frac{1}{2 * |L_{p}|}(\myVec{a}_{p}^{T}\myVec{x}_{i, p} + \omega_{p} - s_{i,p})^{2} \notag \\
    &\qquad \qquad + \sum_{j \in \mySet{N}_p} \myVec{\mu}_{p}^{T}(\myVec{a}_{p} - \myVec{a}_{j}) + \frac{\rho_{p}}{2}\|\myVec{y}_{p} - \myVec{y}_{j}\|_2^2 \notag \\
    &\qquad \qquad + \sum_{j \in \mySet{N}_p} \myVec{\lambda}_{p}^{T}(\myVec{a}_{p} - \myVec{a}_{j}) + \frac{\beta_{p}}{2}\|\myVec{y}_{p} - \myVec{y}_{j}\|_2^2,
    \label{eqn:LinearRegAugLag}
\end{align}
where $\rho > 0$, $\beta > 0$ are a fixed hyperparameters and $\myVec{\mu}_{p}$ $\myVec{\lambda}_{p}$ are the dual variables of node $p$. 


Accordingly, \ac{dadmm} (Algorithm~\ref{alg:Algo1}) can be applied to tackling the \ac{dlr} problem in \eqref{eqn:LinearRegObj} with \eqref{eqn:PrimalStep} becoming
\begin{align}
    \myVec{a}_{p}^{(k+1)} &= \myVec{a}_{p}^{(k)} -\alpha_{p}^{(k)} \nabla_{\myVec{a}_{p}} \mathcal{L}_p(\myVec{y}_{p}^{(k)}, \{\myVec{y}_{j,p}\}_{j\in\mySet{N}_p}, \myVec{\lambda}_p^{(k)}) \notag \\
    &=\myVec{a}_{p}^{(k)} -\alpha_{p}^{(k)} \Big(\frac{1}{|L_{p}|} \sum_{\myVec{x}_{i}, s_{i} \in \mySet{D}_p} \myVec{x}_{i,p}\myVec{x}_{i,p}^{T}\myVec{a}_{p}^{(k)} + \myVec{x}_{i,p}\omega_{p}^{(k)} \notag \\ 
    &- \myVec{x}_{i,p}s_{i,p} +  \sum_{j \in \mySet{N}_{p}} \myVec{\mu}_p^{(k)} + \rho_{p}^{(k)}\big(\myVec{a}_{p}^{(k)} - \myVec{a}_{j,p} \big) 
    \Big),
    \label{eqn:LASSOPrimalStep-a}
\end{align}
and
\begin{align}
    \omega_{p}^{(k+1)} &= \omega_{p}^{(k)} -\delta{p}^{(k)} \nabla_{\omega_{p}} \mathcal{L}_p(\myVec{y}_{p}^{(k)}, \{\myVec{y}_{j,p}\}_{j\in\mySet{N}_p}, \myVec{\lambda}_p^{(k)}) \notag \\
    &=\omega_{p}^{(k)} -\delta_{p}^{(k)} \Big(\frac{1}{|\mySet{D}_{p}|} \sum_{\myVec{x}_{i}, s_{i} \in \mySet{D}_p} \myVec{a}_{p}^{(k)T}\myVec{x}_{i,p} + \omega_{p}^{(k)} + s_{i,p}\notag \\ 
    &\qquad+  \sum_{j \in \mySet{N}_{p}} \myVec{\lambda}_p^{(k)} + \beta_{p}^{(k)}\big(\omega_{p}^{(k)} - \omega_{j,p} \big) 
    \Big),
    \label{eqn:LASSOPrimalStep-omega}
\end{align}
by each agent $p$, with $\alpha > 0$ and $\delta > 0$ being step-sizes.
The agent then shares $\myVec{y}_{p}^{(k+1)} = \{\myVec{a}^{(k+1)}, \omega^{(k+1)}\}$ with its neighbours, who update their local copies. The iteration is concluded by having all agents update their dual variables with step-sizes $\eta > 0$ and $\gamma > 0$ via
\begin{align}
    \myVec{\mu}_{p}^{(k+1)} &= \myVec{\mu}_{p}^{(k)} +\eta_{p}^{(k)} \nabla_{\myVec{\mu}_{p}} \mathcal{L}_p(\myVec{y}_{p}^{(k+1)}, \{\myVec{y}_{j,p}\}_{j\in\mySet{N}_p}, \myVec{\mu}_p^{(k)}) \notag \\
    &= \myVec{\mu}_{p}^{(k)} +\eta_{p}^{(k)}\sum_{j \in \mySet{N}_p}\big(\myVec{a}_{p}^{(k+1)} - \myVec{a}_{j,p}\big).
    \label{eqn:LinearDualStep-mu}
\end{align}
and
\begin{align}
    \myVec{\lambda}_{p}^{(k+1)} &= \myVec{\lambda}_{p}^{(k)} +\gamma_{p}^{(k)} \nabla_{\myVec{\lambda}_{p}} \mathcal{L}_p(\myVec{y}_{p}^{(k+1)}, \{\myVec{y}_{j,p}\}_{j\in\mySet{N}_p}, \myVec{\lambda}_p^{(k)}) \notag \\
    &= \myVec{\lambda}_{p}^{(k)} +\gamma_{p}^{(k)}\sum_{j \in \mySet{N}_p}\big(\omega_{p}^{(k+1)} - \omega_{j,p}\big).
    \label{eqn:LinearDualStep-lambda}
\end{align}



\subsubsection{Unfolding \ac{dadmm}}
By  unfolding the above \ac{dadmm} steps following the method described in Section~\ref{sec:Unfolded}, the iterative steps are viewed as layers of a \ac{dnn}. 
In this problem we formulate the set of learnable hyperparameters on two cases, the first case is when the graph is fixed in time and the second case is when all the agents shares the same learnable hyperparameters.
\begin{enumerate}[label={\arabic*)}]
    \item \label{sssec:linear-case1} 
    {\em Agent-Specific Hyperparameters:}
In this case the set of learnable hyperparameters of agent $p$ is $\myVec{\theta}_{p} = \{\alpha_{p}^{k}, \rho_{p}^{k}, \delta_{p}^{k}, \beta_{p}^{k}, \eta_{p}^{k}, \gamma_{p}^{k}\}_{k=1}^{T}$ (resulting in $6\cdot P \cdot T$ trainable parameters).
    \item {\em Shared Hyperparameters:}\label{sssec:linear-case2}
In this case we neglect the dependency of our model on the number of agents $P$. Therefore, the set of learnable hyperparameters $\myVec{\theta} = \{\alpha^{k}, \rho^{k}, \beta^{k}, \eta^{k}, \gamma^{k}\}_{k=1}^{T}$, will be the same for all agents  (i.e., $6 \cdot T$ trainable parameters). 
\end{enumerate} 
The main difference between these two cases is with the dependency of the learnable hyperparameters set on $P$. 
 The resulting trainable architecture (for an agent-specific parameterization) is illustrated in Fig.~\ref{fig:LinearReg-block-diagram}.

\vspace{-0.2cm}
\subsection{Numerical Evaluation}
\label{ssec:linear-results} 
\vspace{-0.15cm}
%

\begin{figure}
    \centering
    \includegraphics[width=\columnwidth]{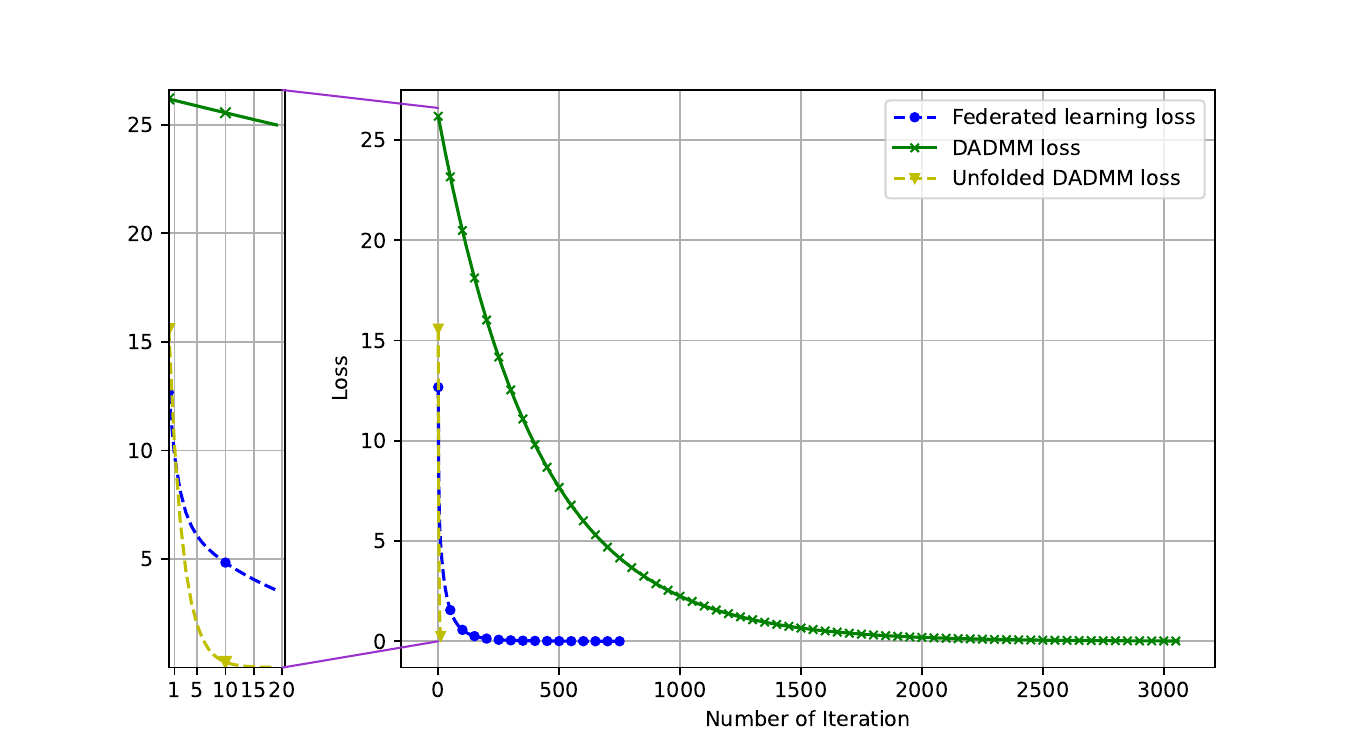}
    \caption{Loss per iteration; $P=5$}
    \label{fig:linear-loss-P=5}
\end{figure}

\begin{figure}
    \centering
    \includegraphics[width=\columnwidth]{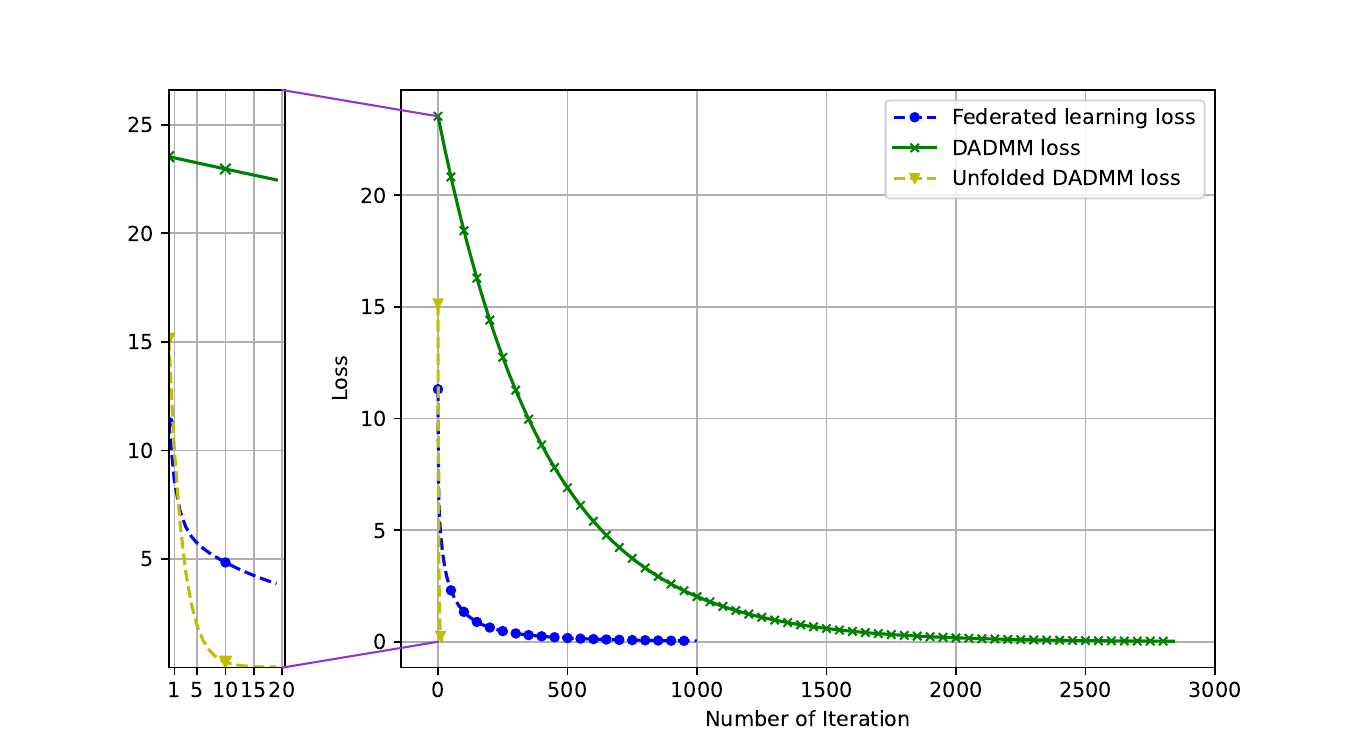}
    \caption{Loss per iteration; $P=12$}
    \label{fig:linear-loss-P=12}
\end{figure}

\begin{figure}
    \centering
    \includegraphics[width=\columnwidth]{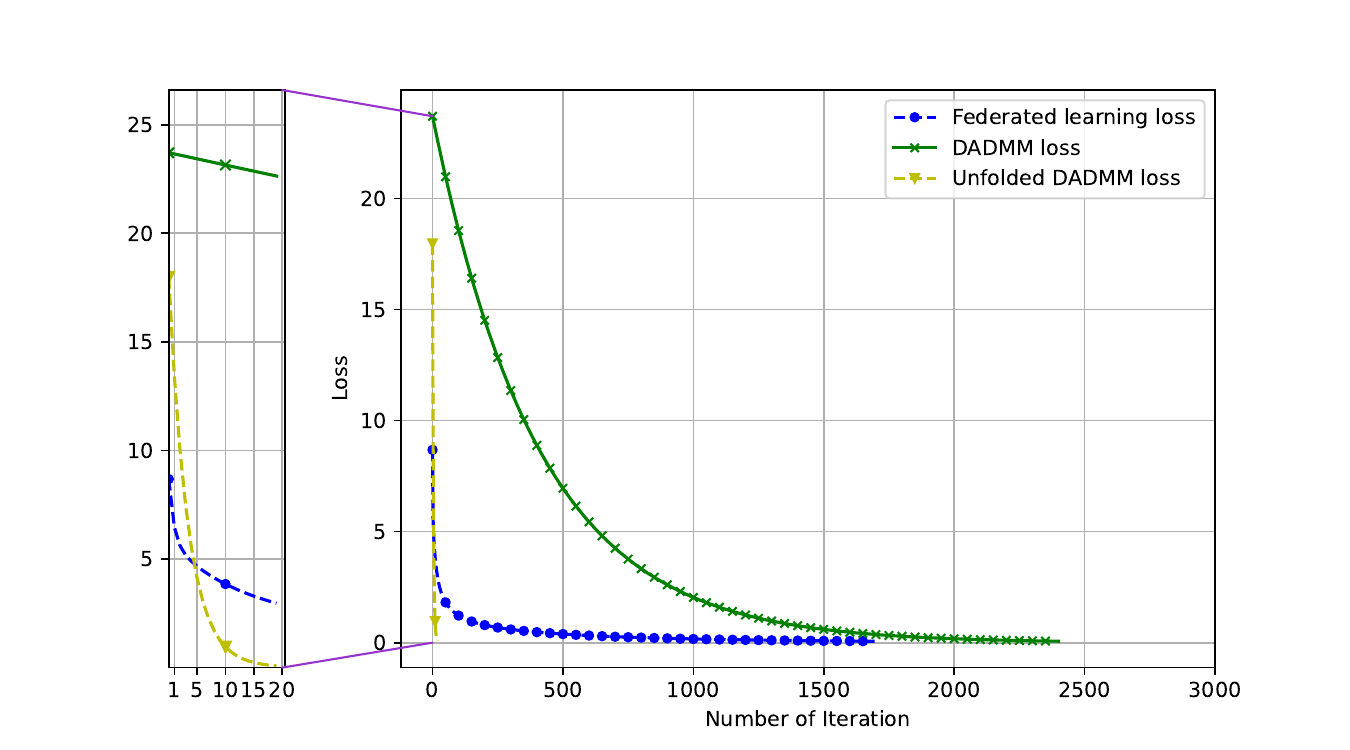}
    \caption{Loss per iteration; $P=20$}
    \label{fig:linear-loss-P=20}
\end{figure}

\begin{figure*}
    \centering
    \minipage{0.48\textwidth}
    \includegraphics[width=\linewidth]{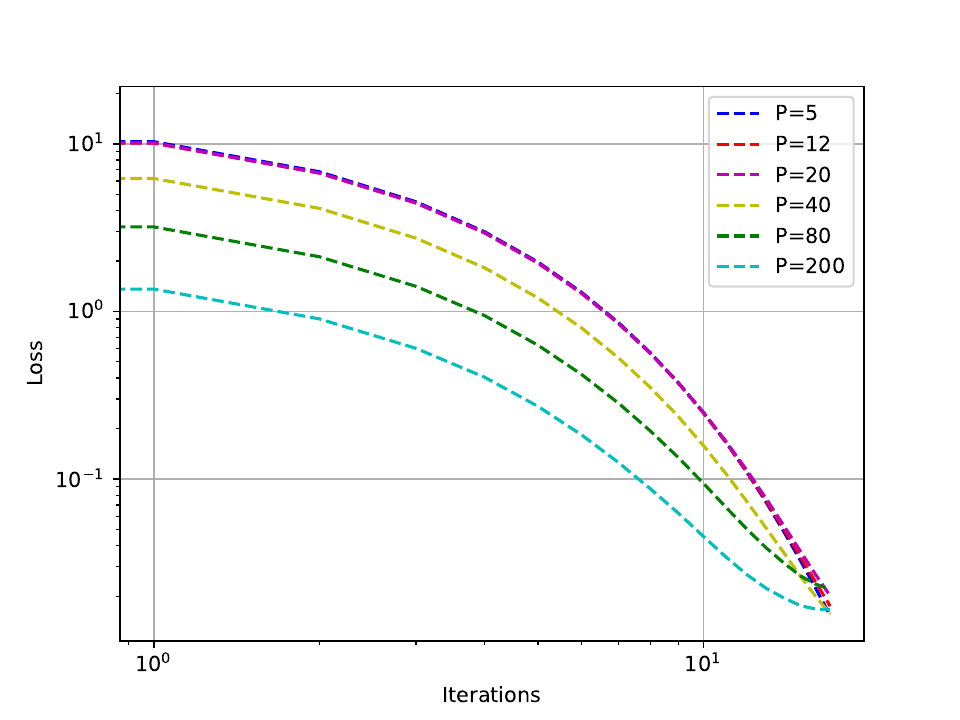}
    \caption{Unfolded \ac{dadmm} with shared hyperparameters loss vs. iteration when trained on graph with $P=5$ agents and tested on bigger graphs $P=\{12, 20, 40, 80, 200\}$.}
    \label{fig:same-hyperparameters-diff-graphs}
    \endminipage\hfill
    \minipage{0.48\textwidth}
    \includegraphics[width=\linewidth]{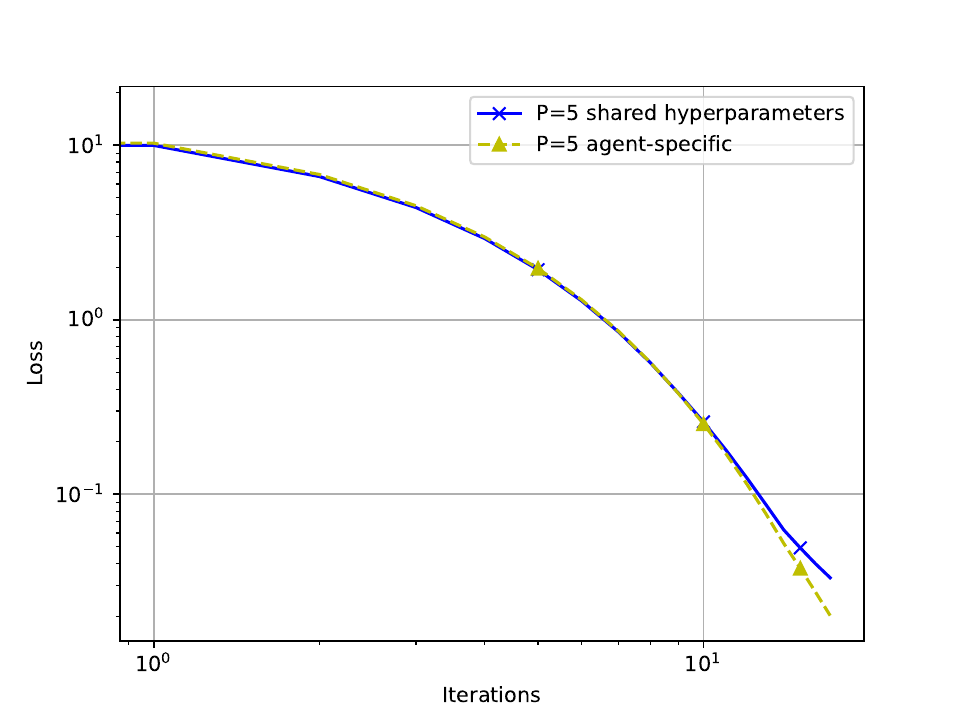}
    \caption{Unfolded \ac{dadmm} loss vs. iteration for agent-specific and shared hyperparameters for graph with $P=5$ agents.}
    \label{fig:same-hyperparameters-compare}
    \endminipage\hfill
\end{figure*}

We numerically evaluate the proposed unfolded D-ADMM algorithm with the two configurations discussed above, i.e.,  $(i)$ Agent-Specific Hyperparameters; $(ii)$ Shared Hyperparameters. 
For both cases we consider the learning of a handwritten digit classifier,  with observations for each agent taken from MNIST dataset as $\myVec{b}_{p} = \{\myVec{x}_{i, p}. s_{i,p}\}_{i=1}^{L_p}$. The desired output $\bar{\myVec{y}}$ is a linear regression model that is applicable for MNIST.

 \subsubsection{Agent-Specific Hyperparameters}
 \label{sssec:linear-case1} 
 Here, we compare unfolded \ac{dadmm} with conventional \ac{dadmm}, where we used fixed hyperparameters manually tuned based on empirical trials to systematically achieve convergence. We also compare distributed optimizers with training the same model using conventional \ac{fl} being a conventional framework for studying distributed machine learning. Note that \ac{fl}  assumes a centralized server that has direct links to each of the $P$ agents. Here, each agent implements $20$ local training iterations using Adam before communicating with the server for synchronization. 

We simulate a communication network using the Erdos-Renyi graph model with proper coloring. The labeled dataset as in \eqref{eqn:dataset} is comprised of $L = L_{p} \cdot P$ training samples, where each agent has local data of size $L_{p} = 200$.
We set the graph size to be $P \in \{5, 12, 20\}$. For each configuration, we apply Algorithm~\ref{alg:L2O} to optimize $\myVec{\theta}$ with $T=20$ iterations using Adam.  
All considered algorithms are evaluated over $200$ test observations.



The results for the settings of $P \in \{5, 12, 20\}$ are depicted in Figs.~\ref{fig:linear-loss-P=5}-\ref{fig:linear-loss-P=20}, respectively. It is observed that the proposed unfolded D-ADMM improves not only upon fixed hyperparameters D-ADMM, but also over \ac{fl}. \textcolor{NewColor}{Specifically, our unfolded \ac{dadmm} achieves improved learned models within limited communication rounds, demonstrating its ability to  improve performance within a pre-defined and limited $T$. This allows it to outperform for small values of $T$ conventional \ac{fl}, that has the advantage of an additional centralized orchestration, yet its learning procedure utilizes fixed hyperparameter patterns for which convergence is achieved for asymptotically large $T$.  }
 Fig.~\ref{fig:linear-loss-P=5} shows a communications reduction by a factor of $\times 154$ and $\times 38$ ($20$ iterations vs. $3080$ and $760$ iterations respectively), while in Fig.~\ref{fig:linear-loss-P=12} the corresponding reduction is by $\times 142$ and $\times 50$ ($20$ iterations vs. $2840$ and $1000$ iterations respectively). In Fig.~\ref{fig:linear-loss-P=20}, we observe a reduction in communication rounds by $\times 120$ and $\times 85$ ($20$ iterations vs. $2400$ and $1700$ iterations respectively). 

 \subsubsection{Shared Hyperparameters}
 \label{sssec:linear-case2}
 Next, we evaluate unfolded \ac{dadmm} with shared hyperparameters. Our main aim here is to evaluate the transferability induced by unfolding with shared hyperparameters, and particularly the ability to train with graph and operate reliably on another graph. 
%
We simulate a communication network using the Erdos-Renyi graph model with proper coloring, and set the graph during training to be comprised of merely $P=5$ nodes. We apply Algorithm~\ref{alg:L2O} to optimize $\myVec{\theta}$ with $T=20$ iterations based on a labeled dataset as in \eqref{eqn:dataset}, where each agent has data with size $L_{p} = 200$, and they all share the same hyperparameters set. 

While we train on a small network, we evaluate the learned hyperparameters set on bigger networks with $P \in \{12, 20, 40, 80, 200 \}$, where all agents has the same learned hyperparameters set. The resulting performance achieved per iteration is depicted in Fig.~\ref{fig:same-hyperparameters-diff-graphs}. It is observed that for the case which all the agents have the same hyperparameters, our proposed method is robust to changes in the number of agents in time, and that one can successfully apply an unfolded \ac{dadmm} trained with small graphs to larger graphs.
  
We next show that the usage of shared hyperparameters does not lead to a notable performance degradation compared to agent-specific hyperparameters, when evaluated and trained on the same graph. To that aim, in Fig.~\ref{fig:same-hyperparameters-compare} we compare the two approaches over the graph used for training, i.e., $P=5$. The figure shows that 
the performance degradation due to reusing hyperparameters, which allows transferability to larger graphs, comes at the cost of only a minor degradation compared with having each agent has its own hyperparameters.

\vspace{-0.2cm}
\section{Conclusions}
\label{ssec:conclusions}
\vspace{-0.15cm} 

In this work, we proposed a data-aided method for rapid and interpretable distributed optimization. Our approach first unfolds the \ac{dadmm} optimizer of each agent to reach consensus using a fixed small number of iterations. Then, we use the data to tune the hyperparameters of each agent at each iteration, which can either be shared between the agents or learned for each agent separately. We specialized our unfolded \ac{dadmm} for distributed sparse recovery and for distributed machine learning, where we showed the notable gains of the proposed methodology in enabling high performance distributed optimization with few communication rounds. 

\bibliographystyle{IEEEtran}
\bibliography{IEEEabrv,refs}

\begin{thebibliography}{10}
\providecommand{\url}[1]{#1}
\csname url@samestyle\endcsname
\providecommand{\newblock}{\relax}
\providecommand{\bibinfo}[2]{#2}
\providecommand{\BIBentrySTDinterwordspacing}{\spaceskip=0pt\relax}
\providecommand{\BIBentryALTinterwordstretchfactor}{4}
\providecommand{\BIBentryALTinterwordspacing}{\spaceskip=\fontdimen2\font plus
\BIBentryALTinterwordstretchfactor\fontdimen3\font minus
  \fontdimen4\font\relax}
\providecommand{\BIBforeignlanguage}[2]{{%
\expandafter\ifx\csname l@#1\endcsname\relax
\typeout{** WARNING: IEEEtran.bst: No hyphenation pattern has been}%
\typeout{** loaded for the language `#1'. Using the pattern for}%
\typeout{** the default language instead.}%
\else
\language=\csname l@#1\endcsname
\fi
#2}}
\providecommand{\BIBdecl}{\relax}
\BIBdecl

\bibitem{noah2023distributed}
Y.~Noah and N.~Shlezinger, ``Distributed {ADMM} with limited communications via
  deep unfolding,'' in \emph{IEEE International Conference on Acoustics, Speech
  and Signal Processing (ICASSP)}, 2023.

\bibitem{nedic2020distributed}
A.~Nedic, ``Distributed gradient methods for convex machine learning problems
  in networks: Distributed optimization,'' \emph{{IEEE} Signal Process. Mag.},
  vol.~37, no.~3, pp. 92--101, 2020.

\bibitem{xin2020general}
R.~Xin, S.~Pu, A.~Nedi{\'c}, and U.~A. Khan, ``A general framework for
  decentralized optimization with first-order methods,'' \emph{Proc. {IEEE}},
  vol. 108, no.~11, pp. 1869--1889, 2020.

\bibitem{yang2019survey}
T.~Yang, X.~Yi, J.~Wu, Y.~Yuan, D.~Wu, Z.~Meng, Y.~Hong, H.~Wang, Z.~Lin, and
  K.~H. Johansson, ``A survey of distributed optimization,'' \emph{Annual
  Reviews in Control}, vol.~47, pp. 278--305, 2019.

\bibitem{boyd2011distributed}
S.~Boyd, N.~Parikh, E.~Chu, B.~Peleato, and J.~Eckstein, ``Distributed
  optimization and statistical learning via the alternating direction method of
  multipliers,'' \emph{Foundations and Trends{\textregistered} in Machine
  learning}, vol.~3, no.~1, pp. 1--122, 2011.

\bibitem{mota2013d}
J.~F. Mota, J.~M. Xavier, P.~M. Aguiar, and M.~P{\"u}schel, ``{D-ADMM}: A
  communication-efficient distributed algorithm for separable optimization,''
  \emph{{IEEE} Trans. Signal Process.}, vol.~61, no.~10, pp. 2718--2723, 2013.

\bibitem{han2012note}
D.~Han and X.~Yuan, ``A note on the alternating direction method of
  multipliers,'' \emph{Journal of Optimization Theory and Applications}, vol.
  155, no.~1, pp. 227--238, 2012.

\bibitem{shi2014linear}
W.~Shi, Q.~Ling, K.~Yuan, G.~Wu, and W.~Yin, ``On the linear convergence of the
  {ADMM} in decentralized consensus optimization,'' \emph{{IEEE} Trans. Signal
  Process.}, vol.~62, no.~7, pp. 1750--1761, 2014.

\bibitem{makhdoumi2017convergence}
A.~Makhdoumi and A.~Ozdaglar, ``Convergence rate of distributed {ADMM} over
  networks,'' \emph{{IEEE} Trans. Autom. Control}, vol.~62, no.~10, pp.
  5082--5095, 2017.

\bibitem{akyildiz2002wireless}
I.~F. Akyildiz, W.~Su, Y.~Sankarasubramaniam, and E.~Cayirci, ``Wireless sensor
  networks: a survey,'' \emph{Computer Networks}, vol.~38, no.~4, pp. 393--422,
  2002.

\bibitem{fischione2011design}
C.~Fischione, P.~Park, P.~D. Marco, and K.~H. Johansson, ``Design principles of
  wireless sensor networks protocols for control applications,'' in
  \emph{Wireless Networking Based Control}.\hskip 1em plus 0.5em minus
  0.4em\relax Springer, 2011, pp. 203--238.

\bibitem{rabbat2004distributed}
M.~Rabbat and R.~Nowak, ``Distributed optimization in sensor networks,'' in
  \emph{International symposium on Information processing in sensor networks},
  2004, pp. 20--27.

\bibitem{shlezinger2022collaborative}
N.~Shlezinger and I.~V. Baji{\'c}, ``Collaborative inference for {AI}-empowered
  {IoT} devices,'' \emph{IEEE Internet of Things Magazine}, vol.~5, no.~4, pp.
  92--98, 2022.

\bibitem{jakovetic2014fast}
D.~Jakoveti{\'c}, J.~Xavier, and J.~M. Moura, ``Fast distributed gradient
  methods,'' \emph{{IEEE} Trans. Autom. Control}, vol.~59, no.~5, pp.
  1131--1146, 2014.

\bibitem{chang2014multi}
T.-H. Chang, M.~Hong, and X.~Wang, ``Multi-agent distributed optimization via
  inexact consensus {ADMM},'' \emph{{IEEE} Trans. Signal Process.}, vol.~63,
  no.~2, pp. 482--497, 2014.

\bibitem{mokhtari2016decentralized}
A.~Mokhtari, W.~Shi, Q.~Ling, and A.~Ribeiro, ``A decentralized second-order
  method with exact linear convergence rate for consensus optimization,''
  \emph{{IEEE} Trans. Signal Inf. Process. Netw.}, vol.~2, no.~4, pp. 507--522,
  2016.

\bibitem{aybat2017distributed}
N.~S. Aybat, Z.~Wang, T.~Lin, and S.~Ma, ``Distributed linearized alternating
  direction method of multipliers for composite convex consensus
  optimization,'' \emph{{IEEE} Trans. Autom. Control}, vol.~63, no.~1, pp.
  5--20, 2017.

\bibitem{lan2020communication}
G.~Lan, S.~Lee, and Y.~Zhou, ``Communication-efficient algorithms for
  decentralized and stochastic optimization,'' \emph{Mathematical Programming},
  vol. 180, no. 1-2, pp. 237--284, 2020.

\bibitem{zhang2017distributed}
G.~Zhang and R.~Heusdens, ``Distributed optimization using the primal-dual
  method of multipliers,'' \emph{{IEEE} Trans. Signal Inf. Process. Netw.},
  vol.~4, no.~1, pp. 173--187, 2017.

\bibitem{sherson2018derivation}
T.~W. Sherson, R.~Heusdens, and W.~B. Kleijn, ``Derivation and analysis of the
  primal-dual method of multipliers based on monotone operator theory,''
  \emph{{IEEE} Trans. Signal Inf. Process. Netw.}, vol.~5, no.~2, pp. 334--347,
  2018.

\bibitem{aysal2008distributed}
T.~C. Aysal, M.~J. Coates, and M.~G. Rabbat, ``Distributed average consensus
  with dithered quantization,'' \emph{{IEEE} Trans. Signal Process.}, vol.~56,
  no.~10, pp. 4905--4918, 2008.

\bibitem{kar2009distributed}
S.~Kar and J.~M. Moura, ``Distributed consensus algorithms in sensor networks:
  Quantized data and random link failures,'' \emph{{IEEE} Trans. Signal
  Process.}, vol.~58, no.~3, pp. 1383--1400, 2009.

\bibitem{zhu2015quantized}
S.~Zhu and B.~Chen, ``Quantized consensus by the {ADMM}: Probabilistic versus
  deterministic quantizers,'' \emph{{IEEE} Trans. Signal Process.}, vol.~64,
  no.~7, pp. 1700--1713, 2015.

\bibitem{shlezinger20221quantized}
N.~Shlezinger, M.~Chen, Y.~C. Eldar, H.~V. Poor, and S.~Cui, ``Quantized
  federated learning,'' \emph{Machine Learning and Wireless Communications},
  pp. 409--432, 2022.

\bibitem{cohen2021serial}
A.~Cohen, N.~Shlezinger, S.~Salamatian, Y.~C. Eldar, and M.~M{\'e}dard,
  ``Serial quantization for sparse time sequences,'' \emph{{IEEE} Trans. Signal
  Process.}, vol.~69, pp. 3299--3314, 2021.

\bibitem{Goodfellow_Deep_Learning}
I.~Goodfellow, Y.~Bengio, and A.~Courville, \emph{Deep Learning}.\hskip 1em
  plus 0.5em minus 0.4em\relax MIT Press, 2016,
  \url{http://www.deeplearningbook.org}.

\bibitem{samikwa2022ares}
E.~Samikwa, A.~Di~Maio, and T.~Braun, ``Ares: Adaptive resource-aware split
  learning for internet of things,'' \emph{Computer Networks}, vol. 218, p.
  109380, 2022.

\bibitem{lee2022radio}
H.-S. Lee, D.-Y. Kim, and J.-W. Lee, ``Radio and energy resource management in
  renewable energy-powered wireless networks with deep reinforcement
  learning,'' \emph{{IEEE} Trans. Wireless Commun.}, vol.~21, no.~7, pp.
  5435--5449, 2022.

\bibitem{alani2022phishnot}
M.~M. Alani and H.~Tawfik, ``Phishnot: A cloud-based machine-learning approach
  to phishing url detection,'' \emph{Computer Networks}, vol. 218, p. 109407,
  2022.

\bibitem{cohen2024sinr}
Y.~Cohen, T.~Gafni, R.~Greenberg, and K.~Cohen, ``{SINR}-aware deep
  reinforcement learning for distributed dynamic channel allocation in
  cognitive interference networks,'' \emph{arXiv preprint arXiv:2402.17773},
  2024.

\bibitem{zhang2023ai}
P.~Zhang, N.~Chen, S.~Shen, S.~Yu, N.~Kumar, and C.-H. Hsu, ``{AI}-enabled
  space-air-ground integrated networks: Management and optimization,''
  \emph{{IEEE} Netw.}, vol.~28, no.~2, pp. 186--192, 2024.

\bibitem{khoshkholghi2022edge}
M.~A. Khoshkholghi and T.~Mahmoodi, ``Edge intelligence for service function
  chain deployment in {NFV}-enabled networks,'' \emph{Computer Networks}, vol.
  219, p. 109451, 2022.

\bibitem{cerquitelli2023machine}
T.~Cerquitelli, M.~Meo, M.~Curado, L.~Skorin-Kapov, and E.~E. Tsiropoulou,
  ``Machine learning empowered computer networks,'' \emph{Computer networks},
  vol. 230, p. 109807, 2023.

\bibitem{song2022networking}
L.~Song, X.~Hu, G.~Zhang, P.~Spachos, K.~N. Plataniotis, and H.~Wu,
  ``Networking systems of {AI}: On the convergence of computing and
  communications,'' \emph{{IEEE} Internet Things J.}, vol.~9, no.~20, pp.
  20\,352--20\,381, 2022.

\bibitem{chen2024learned}
Y.~Chen, L.~Abrahamyan, H.~Sahli, and N.~Deligiannis, ``Learned parameter
  compression for efficient and privacy-preserving federated learning,''
  \emph{IEEE Open Journal of the Communications Society}, vol. 230, pp.
  3503--3516, 2024.

\bibitem{lee2023task}
H.~Lee and S.-W. Kim, ``Task-oriented edge networks: Decentralized learning
  over wireless fronthaul,'' \emph{{IEEE} Internet Things J.}, vol.~11, no.~9,
  pp. 15\,540--15\,556, 2024.

\bibitem{malka2022decentralized}
M.~Malka, E.~Farhan, H.~Morgenstern, and N.~Shlezinger, ``Decentralized
  low-latency collaborative inference via ensembles on the edge,'' \emph{arXiv
  preprint arXiv:2206.03165}, 2022.

\bibitem{vatter2023evolution}
J.~Vatter, R.~Mayer, and H.-A. Jacobsen, ``The evolution of distributed systems
  for graph neural networks and their origin in graph processing and deep
  learning: A survey,'' \emph{ACM Computing Surveys}, vol.~56, no.~1, pp.
  1--37, 2023.

\bibitem{wu2020comprehensive}
Z.~Wu, S.~Pan, F.~Chen, G.~Long, C.~Zhang, and S.~Y. Philip, ``A comprehensive
  survey on graph neural networks,'' \emph{{IEEE} Trans. Neural Netw. Learn.
  Syst.}, vol.~32, no.~1, pp. 4--24, 2020.

\bibitem{shen2020graph}
Y.~Shen, Y.~Shi, J.~Zhang, and K.~B. Letaief, ``Graph neural networks for
  scalable radio resource management: Architecture design and theoretical
  analysis,'' \emph{{IEEE} J. Sel. Areas Commun.}, vol.~39, no.~1, pp.
  101--115, 2020.

\bibitem{rusek2020routenet}
K.~Rusek, J.~Su{\'a}rez-Varela, P.~Almasan, P.~Barlet-Ros, and
  A.~Cabellos-Aparicio, ``Routenet: Leveraging graph neural networks for
  network modeling and optimization in {SDN},'' \emph{{IEEE} J. Sel. Areas
  Commun.}, vol.~38, no.~10, pp. 2260--2270, 2020.

\bibitem{hadou2023stochastic}
S.~Hadou, N.~NaderiAlizadeh, and A.~Ribeiro, ``Stochastic unrolled federated
  learning,'' \emph{arXiv preprint arXiv:2305.15371}, 2023.

\bibitem{lee2019deep}
H.~Lee, S.~H. Lee, and T.~Q. Quek, ``Deep learning for distributed
  optimization: Applications to wireless resource management,'' \emph{{IEEE} J.
  Sel. Areas Commun.}, vol.~37, no.~10, pp. 2251--2266, 2019.

\bibitem{shlezinger2020model}
N.~Shlezinger, J.~Whang, Y.~C. Eldar, and A.~G. Dimakis, ``Model-based deep
  learning,'' \emph{Proc. {IEEE}}, vol. 111, no.~5, pp. 465--499, 2023.

\bibitem{chen2022learning}
T.~Chen, X.~Chen, W.~Chen, Z.~Wang, H.~Heaton, J.~Liu, and W.~Yin, ``Learning
  to optimize: A primer and a benchmark,'' \emph{The Journal of Machine
  Learning Research}, vol.~23, no.~1, pp. 8562--8620, 2022.

\bibitem{shlezinger2023model}
N.~Shlezinger and Y.~C. Eldar, ``Model-based deep learning,'' \emph{Foundations
  and Trends{\textregistered} in Signal Processing}, vol.~17, no.~4, pp.
  291--416, 2023.

\bibitem{shlezinger2022model}
N.~Shlezinger, Y.~C. Eldar, and S.~P. Boyd, ``Model-based deep learning: On the
  intersection of deep learning and optimization,'' \emph{{IEEE} Access},
  vol.~10, pp. 115\,384--115\,398, 2022.

\bibitem{monga2021algorithm}
V.~Monga, Y.~Li, and Y.~C. Eldar, ``Algorithm unrolling: Interpretable,
  efficient deep learning for signal and image processing,'' \emph{{IEEE}
  Signal Process. Mag.}, vol.~38, no.~2, pp. 18--44, 2021.

\bibitem{mogilipalepu2021federated}
K.~K. Mogilipalepu, S.~K. Modukuri, A.~Madapu, and S.~P. Chepuri, ``Federated
  deep unfolding for sparse recovery,'' in \emph{European Signal Processing
  Conference (EUSIPCO)}, 2021, pp. 1950--1954.

\bibitem{nakai2022deep}
A.~Nakai-Kasai and T.~Wadayama, ``Deep unfolding-based weighted averaging for
  federated learning under heterogeneous environments,'' \emph{arXiv preprint
  arXiv:2212.12191}, 2022.

\bibitem{shlezinger2022discriminative}
N.~Shlezinger and T.~Routtenberg, ``Discriminative and generative learning for
  linear estimation of random signals [lecture notes],'' \emph{{IEEE} Signal
  Process. Mag.}, vol.~40, no.~6, pp. 75--82, 2023.

\bibitem{zhu2009distributed}
H.~Zhu, G.~B. Giannakis, and A.~Cano, ``Distributed in-network channel
  decoding,'' \emph{{IEEE} Trans. Signal Process.}, vol.~57, no.~10, pp.
  3970--3983, 2009.

\bibitem{cervino2023training}
J.~Cervino, L.~Ruiz, and A.~Ribeiro, ``Training graph neural networks on
  growing stochastic graphs,'' in \emph{IEEE International Conference on
  Acoustics, Speech and Signal Processing (ICASSP)}, 2023.

\bibitem{shlezinger2020deepsic}
N.~Shlezinger, R.~Fu, and Y.~C. Eldar, ``Deep{SIC}: Deep soft interference
  cancellation for multiuser {MIMO} detection,'' \emph{{IEEE} Trans. Wireless
  Commun.}, vol.~20, no.~2, pp. 1349--1362, 2020.

\bibitem{kingma2014adam}
D.~P. Kingma and J.~Ba, ``Adam: A method for stochastic optimization,''
  \emph{arXiv preprint arXiv:1412.6980}, 2014.

\bibitem{boyd2004convex}
S.~P. Boyd and L.~Vandenberghe, \emph{Convex optimization}.\hskip 1em plus
  0.5em minus 0.4em\relax Cambridge university press, 2004.

\bibitem{shlezinger2022deep}
N.~Shlezinger, A.~Amar, B.~Luijten, R.~J. van Sloun, and Y.~C. Eldar, ``Deep
  task-based analog-to-digital conversion,'' \emph{{IEEE} Trans. Signal
  Process.}, vol.~70, pp. 6021--6034, 2022.

\bibitem{khobahi2021lord}
S.~Khobahi, N.~Shlezinger, M.~Soltanalian, and Y.~C. Eldar, ``Lo{RD-N}et:
  Unfolded deep detection network with low-resolution receivers,'' \emph{{IEEE}
  Trans. Signal Process.}, vol.~69, pp. 5651--5664, 2021.

\bibitem{shlezinger2020communication}
N.~Shlezinger, S.~Rini, and Y.~C. Eldar, ``The communication-aware clustered
  federated learning problem,'' in \emph{IEEE International Symposium on
  Information Theory (ISIT)}, 2020, pp. 2610--2615.

\bibitem{samuel2019learning}
N.~Samuel, T.~Diskin, and A.~Wiesel, ``Learning to detect,'' \emph{{IEEE}
  Trans. Signal Process.}, vol.~67, no.~10, pp. 2554--2564, 2019.

\bibitem{lavi2023learn}
O.~Lavi and N.~Shlezinger, ``Learn to rapidly and robustly optimize hybrid
  precoding,'' \emph{{IEEE} Trans. Commun.}, vol.~71, no.~10, pp. 5814--5830,
  2023.

\bibitem{johnston2021admm}
J.~Johnston, Y.~Li, M.~Lops, and X.~Wang, ``{ADMM-N}et for communication
  interference removal in stepped-frequency radar,'' \emph{{IEEE} Trans. Signal
  Process.}, vol.~69, pp. 2818--2832, 2021.

\bibitem{zhou2021admm}
C.~Zhou and M.~R. Rodrigues, ``An {ADMM} based network for hyperspectral
  unmixing tasks,'' in \emph{IEEE International Conference on Acoustics, Speech
  and Signal Processing (ICASSP)}, 2021, pp. 1870--1874.

\bibitem{wang2022efficient}
M.~Wang, S.~Wei, Z.~Zhou, J.~Shi, and X.~Zhang, ``Efficient {ADMM} framework
  based on functional measurement model for mmw 3-{D SAR} imaging,''
  \emph{{IEEE} Trans. Geosci. Remote Sens.}, vol.~60, pp. 1--17, 2022.

\bibitem{fan2023fast}
W.~Fan, S.~Liu, C.~Li, and Y.~Huang, ``Fast direct localization for millimeter
  wave {MIMO} systems via deep admm unfolding,'' \emph{{IEEE} Wireless Commun.
  Lett.}, vol.~12, no.~4, pp. 748--752, 2023.

\bibitem{shah2024optimization}
S.~B. Shah, P.~Pradhan, W.~Pu, R.~Randhi, M.~R. Rodrigues, and Y.~C. Eldar,
  ``Optimization guarantees of unfolded {ISTA} and {ADMM} networks with smooth
  soft-thresholding,'' \emph{{IEEE} Trans. Signal Process.}, 2024, early
  access.

\bibitem{eldar2012compressed}
Y.~C. Eldar and G.~Kutyniok, \emph{Compressed sensing: theory and
  applications}.\hskip 1em plus 0.5em minus 0.4em\relax Cambridge University
  Press, 2012.

\bibitem{hamilton2017inductive}
W.~Hamilton, Z.~Ying, and J.~Leskovec, ``Inductive representation learning on
  large graphs,'' \emph{Advances in neural information processing systems},
  vol.~30, 2017.

\bibitem{friedlander2012dual}
M.~Friedlander and M.~Saunders, ``A dual active-set quadratic programming
  method for finding sparse least-squares solutions,'' \emph{Online, University
  of British Columbia, BC, Canada}, 2012.

\bibitem{gafni2022federated}
T.~Gafni, N.~Shlezinger, K.~Cohen, Y.~C. Eldar, and H.~V. Poor, ``Federated
  learning: A signal processing perspective,'' \emph{{IEEE} Signal Process.
  Mag.}, vol.~39, no.~3, pp. 14--41, 2022.

\bibitem{ayache2023walk}
G.~Ayache, V.~Dassari, and S.~El~Rouayheb, ``Walk for learning: A random walk
  approach for federated learning from heterogeneous data,'' \emph{{IEEE} J.
  Sel. Areas Commun.}, vol.~41, no.~4, pp. 929--940, 2023.

\end{thebibliography}

\end{document}